\documentclass[a4paper,11pt]{article}
\usepackage{amsfonts}
\usepackage{amssymb}
\usepackage{amsmath}
\usepackage{cite}
\newtheorem{theo}{Theorem}[section]
\newtheorem{lema}[theo]{Lemma}
\newtheorem{rema}[theo]{Remark}
\newtheorem{prop}[theo]{Proposition}
\newtheorem{afir}{Affirmation}[section]

\newtheorem{defi}[theo]{Definition}	
\newcommand{\hypos}{L,M,T,\delta, (\alpha,\beta)}

\newcommand{\pb}{\Gamma_{\alpha}}

\newcommand{\B}{\mathbb{B}}
\newcommand{\Q}{\mathbb{Q}}
\newcommand{\R}{\mathbb{R}}
\newcommand{\D}{\mathbb{D}}

\newcommand{\T}{\mathbb{T}^1}
\newcommand{\cinf}{C^\infty}

\newcommand{\C}{\mathbb{C}}
\newcommand{\N}{\mathbb{N}}
\newcommand{\eps}{\varepsilon}

\newcommand{\tal}{\theta+\alpha}

\newcommand{\Z}{\mathbb{Z}}
\newcommand{\K}{\aleph}
\newcommand{\gab}{\Gamma_{\alpha,\beta}}
\newcommand{\ga}{\Gamma_{\alpha}}

\newcommand{\Pru}{\mathit{Proof. }}

\begin{document}
\title{On the persistence of invariant curves for Fibered Holomorphic Transformations}
\author{Mario PONCE\\  Pontificia Universidad Cat\'{o}lica de Chile}
\maketitle
\begin{abstract}
We consider the problem of the persistence of invariant curves for analytical fibered holomorphic transformations. We define a fibered rotation number associated to an invariant curve. We show that an invariant curve with a prescribed fibered rotation number persists under small perturbations on the dynamics provided that the pair of rotation numbers verifies a Brjuno type arithmetical condition.  
\end{abstract}\


\section{Introduction}

 Skew-product transformations over irrational rotations have been widely studied as a source of  examples of interesting dynamics and for modeling fundamental phenomena in mathematical physics (c.f. quasi-periodic Schr\"odinger cocycle, see \cite{AVKR06, KRIK99}). In the former direction, Furstenberg \cite{FURS61} constructed the first example of a minimal non ergodic conservative diffeomorphism of the torus as a fibered circle diffeomorphism. Latter on, the basis for the study of these maps were stablished by   M.Herman \cite{HERM83}, who defined in particular the notion of the  fibered rotation number (see also \cite{JOMO82}). In the recent years, this theory was relaunched  mainly by the works of G.Keller, J.Stark and T.J\"ager \cite{JAKE06,JAST06}, who established (among other things) a Poincar\'e-like classification relating the fibered rotation number to the existence of invariant graphs. 

In his doctoral thesis \cite{SEST97} (see also \cite{SEST99}), O.Sester studied hyperbolic fibered  polynomials,  successfully generalizing the classical notions of Julia set, Green's function, and the principal cardioid of the Mandelbrot set in the parameter space. Closely related works due to M.Jonsson \cite{JONS99,JONS00} are also important contributions to the subject.\\

In this work we study  fibered holomorphic dynamics.  More precisely, given an irrational number 
$\alpha > 0$, an open simply connected subset $U$ of $\C$, and a real number $\delta>0$, we consider transformations of the form
\begin{eqnarray*} 
F:B_{\delta}\times U &\longrightarrow& B_{\delta}\times \C\\
 (\theta,z)&\longmapsto& \big( \tal,f_{\theta}(z)\big),
 \end{eqnarray*}
where $B_{\delta}$ denotes the strip $\big\{\theta \in \C/\Z\ \big| \ |Im(\theta)|<\delta\big\}$. 
We also assume that this transformation is a holomorphic function as a two variables complex function. In particular, the maps $f_{\theta}:U\to \C$ are holomorphic for every $\theta\in B_{\delta}$. We call $F$ an analytic fibered holomorphic transformation, and we refer to it as a \emph{fht}.

Note that a map $F$ as above has neither fixed nor periodic point. The natural object that plays this  role for the local dynamics is an invariant curve, that is, a holomorhic curve $u:B_{\delta}\to U$  such that, for every $\theta\in B_{\delta}$, 
\begin{equation}\label{introequinv}
 F\big(\theta,u(\theta)\big)=\big(\theta+\alpha,u(\theta+\alpha)\big),
 \end{equation}
which is equivalent to $f_{\theta}\big(u(\theta)\big)\!=\!u(\tal)$. Indeed, the dynamics of $F$ is organized  around this ``fibered fixed point'', thus generalizing the role of a fixed point for the local dynamics of an holomorphic germ $g:(\C,0)\to(\C,0)$ (see \cite{PONC07-1,PONC07t}). We define  an appropriate notion of {\em fibered rotation number}, which appears to be quite useful for studing this dynamics.\\

Fibered holomorphic dynamics may be thought of as  a toy model for understanding the so called Melnikov's problem on the persistence of lower dimensional invariant tori for integrable Hamiltonian systems. Let us recall that, concerning this problem, celebrated results due to H.Eliasson \cite{ELIA88} and J.Bourgain \cite{BOUR97} establish the persistence of lower-dimensional invariant tori provided that the tangencial and normal frequencies verify a diophantine condition $DC({\tau})$.   However, it is conjectured  that condition $DC({\tau})$ is not optimal for this result.

Inspired by this result and the corresponding conjecture, we consider the analogous problem for \emph{fht}'s,  namely the persistence of an invariant curve with prescribed fibered rotation number. The main theorem  of this work states that invariant curves persist under small perturbations under a Brjuno type  arithmetical condition $\B_1$ (which is more general than any condition of type $DC(\tau)$) 
between the rotation on the base $\alpha$ and the fibered rotation number (which in our 
setting corresponds to the tangencial and normal frequencies, respectively).

Although \emph{fht}'s correspond to higher-dimensional systems, due to their skew-product structure 
they have a one-dimensional flavor. In this direction, let us recall that the Brjuno condition is  optimal for the linearization problem of both holomorphic germs and analytic circle diffeomorphisms  close to rotations (see \cite{YOCC95,YOCC99}). In view of this, we conjecture that condition $\B_1$ is optimal for our problem. However, as in the classical case, the use of KAM techniques seems to be inapropriate to prove optimality.

Let us finally point out that the problem of the persistence of invariant curves for $\cinf$ fibered holomorphic maps was settled by the author in \cite{PONC07-2}, where it is proven that the optimal arithmetic condition in that setting is a $DC_1$ condition.

  \paragraph{Acknowledgments. } This work corresponds to the central part of the author's \emph{Th\`ese de Doctorat}, prepared at the \emph{Universit\'e Paris-Sud XI}. The author thanks Jean-Christophe Yoccoz for the suggestion of this subject, his advise and constant support to carry out this work. The author also thanks Jes\'us Zapata for a useful indication, as well as Jan Kiwi, Andr\'es Navas and Godofredo Iommi for several pertinent remarks and corrections. This work was partially supported by CONICYT-Chile and the French government. The article was written during a postdoctoral position of the author at Universidad Cat\'olica de Chile and supported by this institution and the PBCT-Conicyt via the Research Project ADI 17 on Low Dimensional Dynamics.


\section{Arithmetical conditions and linearized equation}\label{APPParitmetico}

Trought this work, $\alpha \!>\! 0$ will be an irrational number (it will correspond to the 
rotation angle on the base for our fibered transformation). For a real number $x$ denote 
by $\|x\|$ the distance to the nearest integer, that is, $\| x \| = \min_{p\in \Z} |x-p|$.


\subsection{Arithmetical conditions for a real number} A complete treatement of what follows in this section may be found in  \cite{LANG66,CASS57,SCHM80}. For each $N\in \N$ we define the \emph{worst divisor} 
of $\alpha$ up to the order $N$ by
\begin{equation*}
\Gamma_{\alpha}(N)=\max_{0< |n|\leq N}\big\|n\alpha\big\|^{-1}.
\end{equation*}  
Worst divisors may be completely described in terms of the denominators of the reduced of the continuous fraction of $\alpha$. Indeed, the sequence $q_0=1, q_1, q_2, \dots$ of these 
denominators has an important property, namely, for every $n \geq 1$, 
\begin{equation*}
q_{n+1}=\min\big\{k>q_n\  \big|\  \|k\alpha\|<\|q_n\alpha\|\big\}.
\end{equation*}
Therefore, if $k\in \N$ is such that $q_k\leq N<q_{k+1}$, then 
$\ga(N)=\big\|q_k\alpha\big\|^{-1}$. 
We define the set
\begin{equation*}
	\B=\big\{\alpha\in \T\setminus \Q\    \big|\   \sum_{N\geq 1}\frac{\log \ga(N)}{N^2}<\infty\big\}.
			\end{equation*}
If  $\alpha$ belongs to $\B$ we say that $\alpha$ verifies the Brjuno (arithmetical) condition 
$\B$ (as defined by A. Brjuno, see \cite{BRJU71}). According to H.R\"{u}ssmann (see 
\cite[\S 8]{RUSS02}), this condition is equivalent to the convergence of the series
\begin{equation*}
\sum_{n\geq 1} \frac{\log q_{n+1}}{q_n}.
\end{equation*}
It is well known that $\B$ has full Lebesgue measure, but is a small set from the topological point of view (Baire's category). The following lemma contains a useful equivalence of Brjuno condition, which we will use through this work. 

\begin{lema}[R\"{u}ssmann\cite{RUSS02}]\label{appplemaalfa}
The number $\alpha$ belongs to $\B$ if and only if
\begin{equation*}
\sum_{n\geq0} \frac{\log \ga(2^n)}{2^n}<\infty.
\end{equation*}
\end{lema}  
$\Pru $ For  every $n\geq 2$ one has 
\begin{equation*}
\sum_{i=2^n}^{2^{n+1}-1}\frac{1}{i}-\frac{1}{i+1}< \sum_{i=2^n}^{2^{n+1}-1}\frac{1}{i^2}< \sum_{i=2^n}^{2^{n+1}-1}\frac{1}{i-1}-\frac{1}{i},
\end{equation*}
and therefore
\begin{equation*}
\frac{1}{2^{n+1}}=\frac{1}{2^n}-\frac{1}{2^{n+1}}<\sum_{i=2^n}^{2^{n+1}-1}\frac{1}{i^2}<\frac{1}{2^n-1}-\frac{1}{2^{n+1}-1}<\frac{1}{2^n}.
\end{equation*}
Since $\ga(\cdot)$ is increasing, this implies that
\begin{equation*}
\frac{1}{2}\frac{\log \ga(2^n)}{2^{n}}\leq \sum_{i=2^n}^{2^{n+1}-1}\frac{\log \ga(i)}{i^2}\leq 2\frac{\ga(2^{n+1})}{2^{n+1}},
\end{equation*}
which allows to conclude $\quad_{\blacksquare}$


\subsection{Arithmetical conditions for a pair of real numbers}

Given a real number $\beta$, we say that the pair $(\alpha,\beta)$ is {\em rational} if there exists 
$k \in \Z$ such that $k\alpha\equiv \beta \mod (1)$.  Whe denote by $\mathbb{T}^2_{\mathcal{I}}$ 
the set of non rational pairs. If $(\alpha,\beta)$ belongs to $\mathbb{T}^2_{\mathcal{I}}$ then, 
for a given $N\in \N$, we define its \emph{worst divisor} up to the order $N$ by
\begin{equation*}
\gab(N)=\max_{0\leq |n|\leq N}\big\|n\alpha-\beta\big\|^{-1}.
\end{equation*}
This represents the quality of the approximations of $\beta$ by multiples of $\alpha$. 
Let $\B_1$ be the set defined by 
\begin{equation*}
\B_1=\Big\{(\alpha,\beta)\in \mathbb{T}^2_{\mathcal{I}}\  \big|\   \alpha\in \B\textrm{ and }\sum_{n\geq1}\frac{\log \gab(n)}{n^2}<\infty\Big\},
\end{equation*} 
and let $\B^{\alpha}_1$ be the set of $\beta$ such that $(\alpha,\beta)$ lies in $\B_1$. 
If $\alpha$ belongs to $\B$ then $\B^{\alpha}_1$ is a set of full Lebesgue measure which 
is meager (that is, equal to a countable union of empty interior closed sets). 
 
\begin{lema}
Assuming that $\alpha$ lies in $\B$, the pair $(\alpha,\beta)$ belongs to $\B_1$ if and only if
\begin{equation*}
\sum_{n\geq 0}\frac{\log \gab(2^n)}{2^n}<\infty.
\end{equation*}
\end{lema}
$\Pru $ Analogous to the proof of Lemma \ref{appplemaalfa}$\quad_{\blacksquare}$


\subsection{Cohomological equation}\label{APPPcohomologique}

For more information on this classical subject see for instance \cite{YOCC92}. Given an 
analytic function $\phi:B_{\delta}\to \C$, we look for an analytical solution 
$\psi:\B_{\delta}\to \C$ to the \emph{cohomological equation}
\begin{equation}\label{appppeqco}
\psi(\tal)-\psi(\theta)=\phi(\theta).
\end{equation}
By integrating both sides of this equality, one readily checks that a necessary condition for 
the existence of a solution to this equation is that the mean value of $\phi$ (with respect to the 
Lebesgue measure in $\T\subset B_{\delta}$) is zero, that is, $\int_{\T}\phi(\theta)d\theta = 0.$
A solution to the equation is not unique, but any two solutions differ by a constant.

\begin{prop} \label{eqcoan}
If $\alpha$ belongs to $\B$ and the mean value of $\phi$ is zero, then there exists an analytic function $\psi\!:B_{\delta}\to \C$ which is a solution to {\em (\ref{appppeqco})}. Moreover, such a solution may be taken so that $\int_{\T}\psi(\theta)d\theta=0$ $\quad_{\blacksquare}$.
\end{prop}

 Although condition $\B$ is not optimal for this proposition, it is optimal for the problem of  linearization of holomorphic germs in the neighborhood of indifferent irrational fixed points  \cite{YOCC95}.
\begin{rema}
An easy homotopy type argument shows that the mean value of $\phi$ on the circle $\T_c$ is  
constant as a function of $c\in (-\delta,\delta)$. If the value of this constant is non zero 
then we may envisage equation (\ref{appppeqco}) replacing $\phi$ by 
$\phi-\int_{\T_c}\phi(\theta)d\theta$ for some (any) $c$.
 \end{rema}
\section{Definitions and normal forms}\label{FORMNORM}
 In order to state properly the main theorem of this work we need to introduce some definitions.
\\

Let $F$ be a \emph{fhd} and $u$ an invariant curve for $F$. We consider the circle $\T$ naturally embedded in $B_{\delta}$. In the remaining part of this work  we will always suppose that $\int_{\T}\log \big|\partial_zf_{\theta}\big(u(\theta)\big)\big|d\theta=0$.  (Notice that, since $F$ is injective, the differential $\partial_zf_{\theta}$ is always non zero.) This condition says that the curve is neither attracting nor repulsive at the infinitesimal level. We say that the curve is indifferent. We also suppose that the  application $\theta\mapsto \partial_{z}f_{\theta}\big(u(\theta)\big)$
 is homotopic  in $\C\setminus \{0\}$ to a constant. 
\begin{defi}\label{deffrn}
We define the \emph{fibered rotation number} by
\begin{equation*}
\varrho_{tr}(u)=\frac{1}{2\pi i }\int_{\T}\log \partial_zf_{\theta}\big(u(\theta)\big)d\theta .
\end{equation*}
This number  represents  the average rotation speed of the dynamics around the invariant curve.  Notice that the $\log$ above is well defined $\mod 2\pi i$, and the number $\varrho_{tr}(u)$ is well defined $\mod 1$. 
\end{defi}
\begin{rema}
Notice that, since all of the circles $\T_c=\big\{\theta\in B_{\delta}\ \big| \ Im(\theta)=c\big\}$ are homotopic between them and $F$ is holomorphic, one has 
\begin{equation*}
\varrho_{tr}(u)=\frac{1}{2\pi i}\int_{\T_c} \log \partial_zf\big(\theta,u(\theta)\big)d\theta
\end{equation*}
for every circle $\T_c$. Thus, the indifferent nature of the curve (resp. the fibered rotation number)  may be detected (resp. computed) on any circle $\T_c$, for $|c|<\delta$.  
\end{rema}

Along this work  we will often deal with  holomorphic changes of coordinates $\tilde{H}$  defined from a tubular neighborhood of the zero section $\{z\equiv 0\}_{B_{\delta}}$ to a tubular neighborhood of the invariant curve, say
\begin{equation*}
(\theta,z)\stackrel{\tilde{H}}{\longmapsto}\big(\theta,h_{\theta}(z)\big).
\end{equation*}
The functions $h_{\theta}$ will be  biholomorphic transformations between two topological discs sending the origin to the curve, that is, $h_{\theta}(0)=u(\theta)$. Moreover, the application $\theta\mapsto h_{\theta}(0)$ will be homotopic to a constant in $\C\setminus \{0\}$ (in other words, it will have zero topological degree). Therefore, when conjugating our transformation $F$ by such a $H$ we will get a new transformation $\tilde{F}=H^{-1}\circ F\circ H$ having the zero section $\tilde{u}=\{z\equiv 0\}_{B_{\delta}}$ as an indifferent invariant  curve and  having the same fibered rotation number, that is, $\varrho_{tr}(\tilde{u})=\varrho_{tr}(u)$. In particular, the indifferent nature of the curve  and the fibered rotation number are invariant by this type of conjugacy.  
\\

 Assume $F:B_{\delta}\times U\to B_{\delta}\times \C$ is a \emph{fhd}, $\alpha$ is a Brjuno number, and $u$ is  an indifferent invariant curve for $F$. For $\beta=\varrho_{tr}(u)$ we will solve the (cohomological) equation
 \begin{equation*}
 v(\theta)-v(\theta+\alpha)=2\pi i \beta-\log \partial_zf_{\theta}\big(u(\theta)\big)
 \end{equation*} 
on the strip $B_{\delta}$ (see Proposition \ref{eqcoan}).  Then letting $u_1=\exp (v)$ and $H(\theta,z)=\big(\theta, u(\theta)+u_1(\theta)z\big)$, and performing the corresponding change of coordinates, one gets the following \emph{normal form} for $F$:
 \begin{equation*}
N_F (z,\theta)= H^{-1}\circ F\circ H(\theta,z)=\big(\tal, e^{2\pi i \beta}z+ \rho(\theta,z)\big).
 \end{equation*}
Here $\rho$ is an analytical function vanishing up to the second order at $z=0$. This function $\rho$, and so the \emph{fhd} $N_F$, is defined on the strip $B_{\delta}$ for the $\theta$ variable, and $\rho(\theta,\cdot)$ is defined on a non-constant open set $U_{\theta}=H^{-1}(\{\theta\}\times U)$. \\

 Assume now that $F$ is a \emph{fhd} that may be putted in a normal form 
 \begin{equation*}
 H^{-1}\circ F\circ H^{-1}(\theta,z)=\big(\tal, e^{2\pi i \beta}z+\rho(\theta,z)\big)
 \end{equation*}
by a change of coordinates of the form 
$$H(\theta,z)=\big(\theta,u_0(\theta)+u_1(\theta)z\big),$$
where  $u_0:B_{\delta'}\to U$, $u_1:B_{\delta'}\to\C$ are  analytical functions  for some $\delta'\in (0,\delta]$, and  $u_1$ has zero topological degree.   
  In this case one can readily check that $u_0$ is an indifferent invariant curve, and its fibered rotation number equals $\varrho(u_0)=\beta$.  
 \\
 
Therefore, the existence of  an indifferent invariant curve with transversal rotation number $\beta$ is equivalent, under the Brjuno condition over $\alpha$, to the existence of an analytical fibered affine change of coordinates that puts $F$ in an apropriate normal form. 
 \section{Statement of the problem and result}
Let $F$ be  a \emph{fhd} with an indifferent invariant  curve having a given fibered rotation number, say $\beta\in \R$. As we claimed in the previous Section, if $\alpha \in \B$ then $F$ can be written in a neighborhood of the curve in a normal form $\big(\tal, e^{2\pi i \beta}z+\rho(\theta, z)\big)$, where $\rho$ is an analytic function. The function $\rho(\theta,\cdot)$ vanishes up to the order $2$ at $z=0$, and is convergent for $|z|<r$ for some positive radius $r$.  A small perturbation of such a transformation is defined as being a \emph{fhd} of the form
 \begin{equation*}
 (\theta,z)\longmapsto \big(\tal, \rho_0(\theta)+\big(e^{2\pi i \beta}+\rho_1(\theta)\big)z+\tilde{\rho}(\theta,z)\big),
 \end{equation*}
where $\rho_0,\rho_1$ are small analytic functions, where $\tilde{\rho}$ is an analytic function whose em size is comparable to that of $\rho$, and where $\tilde{\rho}(\theta,\cdot)$ vanishes up to the order $2$ at $z=0$ and is convergent for $|z|<r$.
 \paragraph{Notation for sizes. }
Every function that will appear in this work is a bounded holomorphic function with respect to its variables (defined in some complex open set). Let $f:U\subset \C^k\to \C$ be a bounded holomorphic function from the open set  $U\subset\C^k$, $k\in \{1,2,3\}$, to $\C$. We define the norm $\|f\|_U$ of $f$ as 
 \begin{equation*}
 \|f\|_U=\sup_{z\in U}|f(z)|.
 \end{equation*}
 In general,  $U$  can be thought as a product $B_{\delta}\times D^Z\times D^t$, where  $D^Z,D^t$ are discs in the complex plane. We then may define the oscillation of $f$ with respect to the $z$ variable as
 \begin{equation*}
 osc\big(f\big)_{B_{\delta},D^Z,D^t}=\sup_{\begin{subarray}{1} (\theta,t)\in B_{\delta}\times D^t\\ z_1, z_2 \in D^Z  \end{subarray}}|f(\theta,z_1,t)-f(\theta,z_2,t)|.
 \end{equation*}
  Note that if $z=0$ belongs to  $D^Z$ and $f$ vanishes at this point, then its oscillation is a bound for its norm, that is, 
   \begin{equation*}
\|f\|_{B_{\delta},D_Z,D_t}\leq osc\big(f\big)_{B_{\delta},D_Z,D_t}\leq 2\|f\|_{B_{\delta},D^Z,D^t}  .
\end{equation*}
We will currently omit the corresponding domains of definition when this does not lead to confussion. 

\subsection{$1-$parameter families}

In general KAM results, a perturbation over  an \emph{elliptical dynamics} (rotations of the circle, completely integrable Hamiltonians, etc.) leads to a perturbation of the frequencies which control the dynamics. Thus, in the perturbed situation,  we can not expect to retrieve the same dynamical properties than in the unperturbed elliptic situation. In order to retrieve these original properties, one introduces a $1-$parameter correction, where the parameter is a vector value of the same nature (and dimension) as the related frequencies. In this way, one usually shows that some dynamical properties (linearization, existence of invariant curves with a given rotation number, etc.) are \emph{persistent} in codimension $1$ (see for instance \cite{BOST86}).  In our work, we will only perturb the holomorphic part of the transformations.  Therefore, our frequency corresponds to the fibered rotation number, that is, a complex number.
  \\

Let $F$ be a \emph{fhd} and $u$ an indifferent invariant curve, with $\varrho_{tr}(u)=\beta\in \R$. Let $\Sigma$ be an open set in $\C$.  A \emph{transversal small perturbation} of $F$ is  a 1-complex parameter family $\{F_t\}_{t\in \Sigma}$ of \emph{fhd}'s (i.e., a complex curve in the space of \emph{fhd}'s)  verifying the following properties:  
  
\noindent -- every element $F_t$ is a small perturbation of $F$, 

\noindent -- the  fibered rotation number (even if an invariant curve does not exist) changes along the complex curve $\{F_t\}_{t\in \Sigma}$ (see Theorem \ref{teo} for more explanations). 

We say that the curve $u$ is \emph{persistent} if for any transversal small perturbation $\{F_t\}_{t\in \Sigma}$ there exists a parameter $t^*\in \Sigma$ such that $F_{t^*}$ has an invariant indifferent curve $u^*$ satisfying $\varrho_{tr}(u^*)=\beta$. Roughly speaking, our result says that, except for a complex correction,  any small perturbation of $F$ has an indifferent invariant curve having $\beta$ as its fibered rotation number. However, an arithmetical condition over the pair $(\alpha, \beta)$ is required.
  \\
  
   We say that a  1-parameter family $F_t$ of \emph{fhd}'s is \emph{analytic} if the application 
   \begin{eqnarray*}
   \Sigma\times B_{\delta}\times D(0,r)&\longrightarrow&\C\\
   (t,\theta,z)&\longmapsto& F_t(\theta,z)
   \end{eqnarray*}    
 is holomorphic. We consider the usual development for the elements of this family, namely, 
 \begin{equation}\label{notacionfamilia}
 F_t(\theta,z)=\big(\tal, \rho_{0,t}(\theta)+\big(e^{2\pi i\beta}+\rho_1(\theta)\big)z+\rho_t(\theta,z)\big).
\end{equation}
  To fix ideas, we will suppose that $\rho_t(\theta,\cdot)$ is convergent in the unit disc $\D$. With this notation we can finally state our main result.  
  	\begin{theo}\label{teo}
	For every pair $(\alpha,\beta)$ verifying the  Brjuno condition $\B_{1}$, and for every positive constants $L>1,M,T$ and $\delta$, there exists a real number $\eps^*(L,M,T,\delta,(\alpha,\beta))>0$ and a positive constant  $K_R(\hypos)$ such that the following holds:  if  $\eps$ belongs to  $(0,\eps^*]$ and for an analytic family  $\{F_t\}_{t\in\Sigma}$ there exists a disc  $D(t_0,K_R\eps)\subset\Sigma$ such that  
	\begin{equation}\label{transversality}
		L>\Big|\partial_t\int_{\T}\rho_{1,t}(\theta)d\theta\big|_{t=t_0}\Big|>L^{-1}
		\end{equation}
		and such that, for every $t$ in $D(t_0,K_R\eps)$, 
		\begin{displaymath}
		\left \{ \begin{array}{l}
		\|\rho_{0,t}\|_{B_{\delta}}\leq \eps\\
	\|\rho_{1,t}\|_{B_{\delta}}\leq\eps\\
		\|\partial_z^2\rho_t\|_{B_{\delta},\D}\leq M\\
		\|\partial_t\partial_z\rho_t\|_{B_{\delta},\D}+\|\partial_t^2\rho_{1,t}\|_{B_{\delta}}\leq T, 
		\end{array}\right.
\end{displaymath}
then there exists a parameter  $\bar{t}$ in $D(t_0,K_R\eps)$ such that $F_{\overline{t}}$ has an indifferent invariant curve $u$ which is analytic on the strip  $B_{\frac{\delta}{2}}$. Moreover,  the fibered rotation number of this curve equals $\beta$. Furthermore, the size of  $u$ goes to $0$ when $\eps$ goes to $0$.
	\end{theo}
	Note that the transversality condition is given by (\ref{transversality}). In \cite{PONC07-2} we show an analogous result in the case of $\cinf$ transformations. Indeed, we show that a diophantine condition over the pair $(\alpha,\beta)$ is sufficient in order to get persistence of the invariant curve. We also show that this arithmetical condition is optimal. For a non-diophantine pair $(\alpha,\beta)$ we present a smooth \emph{fhd} having a non-persistent indifferent invariant curve with $\beta$ as fibered rotation number. We conjecture that in the analytic setting, the Brjuno $\B_1$ arithmetical condition is optimal in this sense.  
\paragraph{Plan of the proof. }
The proof of Theorem \ref{teo} will be obtained as an application of the successive conjugacy   method (Newton's algorithm). We will find  an analytic  fibered affine coordinates change that will put one  of the \emph{fhd}'s $F_t$ in the family  in the adequate normal form. As usual in this type of proofs, we must keep the control on several quntities at the same time and perform many different operations. Thus, the whole proof is divided in various parts that we describe now. 
In Section \ref{troncature} we make a review of some simple but useful results in complex analysis that will be used in the demonstration of Theorem \ref{teo}. Some conventions about the names of the constants arising in the work are done in that Section. In Section \ref{lemas} we show two lemmas which represents the core of the proof. These lemmas define the convergence speeds of the iterative process.  Moreover,  they give the size of the width strip loses that we must introduce at each stage. These lemmas establish the relations between the arithmetical properties of the rotations numbers and the characteristics of the iterative process  to be applied, and therefore, between arithmetics and dynamics. In Section \ref{proceso}  we introduce the actual iterative process. In \ref{aleph} we define a constant $\K$ which determines the values given by the lemmas in Section \ref{lemas}. In \ref{siete} we make a detailed description of the algorithm. We introduce and describe the different operations used at each stage.  We also list the order in which these operations will be applied.  In \ref{depapropab} we make some preliminary estimates  and we assume two additional hypothesis in order to reach the  adequate situation which enables us to  start with the  iterative process.   In \ref{realisacion},\ref{estimacionent}, \ref{coup} we perform in a detailed way all the operations and estimates  of a single stage  of the algorithm. We show indeed that these operations allow  us  to get the good estimates  in order to iterate the process and thus,  to approximate the normal form better and better. This is the more delicate part of the proof, since we must  control many quantities that arise as the result of the different operations.  In \ref{convergencia}    we show that our algorithm converges, and thus,  it gives  the desired normal form for one of the  elements $F_t$ in the family.  In \ref{remarcaepsilon} we exhibit a modification of the first stage of the algorithm, in order to obtain more precise results when $\eps$ is very small. In \ref{previous} we remove  the additional hypothesis introduced in  Section \ref{depapropab}. We  show that we can get the desired situation as a consequence of the original hypothesis of Theorem \ref{teo} by means of a previous preparative procedure.   In \ref{teopara} we prove a parametrized version of the main theorem. We expect that this version of the theorem should  be useful in many other applications. 
	\section{Some results on analysis}\label{troncature}
	\subsection{Constants}\label{constantes}
	In this work,  several positive constants with different origin will appear. Some of them have a fixed numerical value. These constants are called \emph{universal} and are denoted by  $C_j$,  where the index represents the order of apparition of the constant in the work. If we denote  a constant by  $C_j$ then this means that this constant is universal. There exists constants which make part of  the statement of the theorem, these constants are $\hypos$.  In a third group, there are constants whose values depend on the values of the constants $\hypos$, and these are denoted by $K_j$. The index represents the order of apparition or the feature of the constant. The special notation $\K$ is reserved to a particular constant belonging to this last group. At some extent, this constant determines completely the proof.
 \subsection{Estimations for truncate analytical functions }\label{troncado}
 Let $f$ be an analytical function defined on the strip  $B_{\delta}$, and let $N$ be a natural number. We define respectively the truncation  up to the order   $N$  and the rest up to the order   $N$ of $f$ by 
 \begin{eqnarray*}
 f|_N(\theta)&=&\sum_{|n|\leq N}\hat{f}(n)e^{2\pi in\theta}\\
 f|_{\overline{N}}(\theta)&=&\sum_{|n|>N}\hat{f}(n)e^{2\pi in\theta},
 \end{eqnarray*}
 where $\hat{f}(n)$ is the $n-$th coefficient of the Fourier series of $f$. 
 These two functions are analytical in the strip  $B_{\delta}$. If  $\delta'\in (0,\delta)$ and  $\delta-\delta'$  is small then there exists a constant  $C_1$ such that 
 \begin{eqnarray}
 |\hat{f}(n)|&\leq&\|f\|_{B_{\delta}}e^{-2\pi |n|\delta}\\\label{estimarestos}
 \|f|_{\overline{N}}\|_{B_{\delta'}}&\leq&C_1\|f\|_{B_{\delta}}\frac{e^{-2\pi N(\delta-\delta')}}{\delta-\delta'}.
 \end{eqnarray}
 The last inequality says that a small loss in the width  of the definition strip gives an estimate on the norm of the rest of $f$ (Cauchy estimate).  We recall here the following simple fact.
 \begin{afir} \label{afirmac} If $f$ is an analytical function whose differential $\partial_{\theta}f$ is bounded in  $B_{\delta}$ then
 \begin{equation}
 \|f\|_{B_{\delta}}\leq \bigg|\int_{\T}f(\theta)d\theta\bigg|+\|\partial_{\theta}f\|_{B_{\delta}}.
 \end{equation}
 \end{afir}
 \subsection{Estimates for the solutions of some truncate cohomological equations}\label{cohomo}
  We consider the following cohomological equations on  $g_1$ and $g_2$
  \begin{eqnarray}
  g_1(\theta)-g_1(\theta+\alpha)&=&p_1(\theta)\\
  e^{2\pi i \beta}g_2(\theta)-g_2(\theta+\alpha)&=&p_2(\theta),
  \end{eqnarray}
  where $\alpha,\beta$ are real numbers,  $p_1,p_2$ are trigonometric polynomials of orders smaller than $N+1$, and $\int_{\T}p_1(\theta)d\theta=0$. Using the Fourier coefficients method we obtain the expressions for the solutions (as explained in the Introduction, for the first equation we just retain the solution with zero mean):
  \begin{eqnarray*}
  g_1(\theta)&=& \sum_{0<|n|\leq N} \frac{\hat{p}_1(n)e^{2\pi i n\theta}}{1-e^{2\pi i n \alpha}}\\
  g_2(\theta)&=&\sum_{|n|\leq N} \frac{\hat{p}_2(n)e^{2\pi i n\theta}}{e^{2\pi i\beta}-e^{-2\pi i n\alpha}}.
  \end{eqnarray*}  
Let $f_1,f_2$ be analytic functions defined on $B_{\delta}$, with $\int_{\T}f_1=0$.  For j... we put  $p_j=f_j\big|_N$, for $j\in \{1,2\}$.  Denoting $\varphi_1(N)=\ga (N), \varphi_2(N)=\gab (N)$ (see Section \ref{APPParitmetico}), we get a constant $C_2$ such that
   \begin{eqnarray*}
   \|g_j\|_{B_{\delta}}&\leq&C_2N\|f_j\|_{B_{\delta}}\varphi_j(N)\\
   \|\partial_{\theta}g_j\|_{B_{\delta}}&\leq&C_2N^2\|f_j\|_{B_{\delta}}\varphi_j(N)
   \end{eqnarray*}
   for $j\in \{1,2\}$. Note that the factor $C_2N$ is far from being optimal, but it will be enough for our purposes.
     \section{Two technical lemmas}\label{lemas}
     In this Section we show two technical lemmas which will provide the definition of two sequences of positive numbers 
$\{l_n\},\{w_n\}$ converging  to $0$ in a very particular way. These sequences will be used as a measure of the rate of convergence of our procedure. These two lemmas will also provide two summable sequences of positive numbers  $\{d_n^0\},\{d_n^1\}$. These sequences  will measure loses in the width  of the analyticity strip in order to get some Cauchy estimates. This Section is the most important  part from the point of view of the relations between dynamics and arithmetic.
Let $\K>2$ be a real parameter and  $\{A_n\}_{n\geq 0},\{B_n\}_{n\geq 0}$ be the non-decreasing sequences defined by  $A_n=\gab(2^n)$ and $B_n=\ga(2^n)$. Under the Brjuno $\B_{1}$ condition these sequences verify
\[
\sum_n\frac{\log A_n}{2^n}<\infty, \qquad \sum_n \frac{\log B_n}{2^n}<\infty .
\] 
\begin{lema}\label{ele1}
For
\begin{equation}
l_n=\frac{e^{-n^2}}{16\K A_nB_n2^n}
\end{equation}
the following holds: 
\begin{enumerate}
\item $\lim_{n\to \infty} l_n=0$.
\item $\sum_n \K 2^nB_nl_n<\frac{1}{8}$.
\item $\K 2^nA_nl_n<\frac{1}{16}$.
\item The relation 
\begin{equation}\label{relad1}
l_{n+1}=\frac{\K l_ne^{-2\pi 2^n d_n^1}}{d_n^1}
\end{equation}
defines a sequence $\{d_n^1\}_{n\geq 0}$ of real positive numbers such that 
\begin{equation}
\sum_n d_n^1 <K_1\big(1+\log \K\big)
\end{equation}
for a constant  $K_1$ depending only on   $\{A_n\},\{B_n\}$.
\end{enumerate}
\end{lema} 
$\Pru$ The proof of the statements  \emph{1., 2., 3.} are straightforward. For \emph{4.} we put   
\begin{eqnarray*}
\hat{\sigma}_n&=&-\log(\K 2^nl_n)\\
&=&\log A_n+\log B_n+\log 16+ n^2,
\end{eqnarray*}
and we define a sequence  $\{\tilde{d}_n^1\}_{n\geq 0}$ by
\begin{equation*}
l_{n+1}=\K l_ne^{-2\pi 2^n \tilde{d}_n^1}.
\end{equation*}
 We obtain that 
\begin{equation*}
2\pi 2^n\tilde{d}_n^1=\hat{\sigma}_{n+1}-\hat{\sigma}_n+\log2\K .
\end{equation*}
As $\hat{\sigma}_n$ is increasing, the numbers $\tilde{d}_n^1$ are positive. Since the  series  $\sum \frac{\hat{\sigma}_n}{2^n}$ is convergent, the series $\sum \tilde{d}_n^1$ is also convergent, and its value is bounded by  $K_1(1+\log \K)$, where  $K_1$ only depends on  $\{A_n\}, \{B_n\}$.

We construct now a sequence  $\{x^1_n\}_{n\geq 0}$ in such a way that the perturbed numbers $d_n^1=\tilde{d}_n^1+x^1_n$ verify point $4.$. For this, we define $x^1_n$ as being the solution of the equation
\begin{equation}\label{luna}
\frac{e^{-2\pi 2^n x^1_n}}{\tilde{d}^1_n+x^1_n}=1.
\end{equation}
It is easy to see that  (\ref{luna}) has a unique solution  $x^1_n$. Moreover, we have  $x_n<s_n$, where  $s_n$ is the solution to  
\begin{equation*}
e^{-2\pi 2^n s_n}=s_n.
\end{equation*}
This equality easily implies that the series $\sum s_n$ converges, thus concluding the proof $\quad_{\blacksquare}$
\begin{lema}\label{we1}
For
\begin{equation}
w_n=\frac{l_n}{\K 4^nA_n},
\end{equation}
where the sequence $\{l_n\}_{n\geq 0}$ is defined as in Lemma \ref{ele1}, the following holds:
\begin{enumerate}
\item $\K w_n(2^nA_n)^2<\frac{1}{32}$.
\item $\sum_n\K w_n2^nA_n<\frac{1}{16}$.
\item The relation
\begin{equation}\label{relad0}
\frac{w_{n+1}}{4}=\frac{w_n\K e^{-2\pi 2^nd_n^0}}{d_n^0}
\end{equation}
defines a  sequence  $\{d_n^0\}_{n\geq 0}$ of positive numbers  such that  
\begin{equation}
\sum_n d_n^0\leq K_2(1+\log \K),
\end{equation}
for a constant  $K_2$ depending  only on $\{A_n\},\{ B_n\}$.
\end{enumerate}
\end{lema}
$\Pru$ The proof of the statements   \emph{1., 2.} are straightforward. Let us show \emph{3.}. As in the proof of Lemma \ref{ele1} we define a sequence  $\{\tilde{d}_n^0\}_{n\geq 0}$ by the relation
\begin{equation*}
\frac{w_{n+1}}{4}=w_n\K e^{-2\pi 2^n\tilde{d}_n^0},
\end{equation*}
that is, 
\begin{equation*}
\tilde{d}^n_0=\frac{1}{2\pi 2^n}\big \{2n+1+\log (A_{n+1}^2B_{n+1})-\log(A_n^2B_n)+\log 32\K\big \}.
\end{equation*}
Therefore, we have a positive and summable sequence. As in Lemma  \ref{ele1}, we can find a real sequence  $\{x_n^0\}_{n\geq 0}$ such that the perturbed numbers   $d_n^0=\tilde{d}_n^0+x_n^0$ verify \emph{3.} $\quad_{\blacksquare}$
\section{An iterative process}\label{proceso}
\subsection{The choice of $\K$}\label{aleph}
The sequences defined in the previous Section depend on the value of the constant $\K$.
 In order to use these sequences in our process we need to choose the value $\K$  being larger than many quantities which depend only on the values of constants of type $K$ (cf. Section \ref{constantes}). This choice can be made a priori  since the type $K$ constants depend only on the values of the constants $\hypos$ that appear in  the statement of Theorem \ref{teo}. However, along the process we will write explicitly the bounds assumed for $\aleph$.   Fixing the value of $\K$, we obtain   two positive sequences $\{d_n^0\}_{n \geq 0},\{d_n^1\}_{n\geq 0}$. We will use these sequences as loses in the width of the analyticity strip. These sequences verify  
\begin{equation*}
\sum_{n\geq 0}\max\big(d_n^0,d_n^1\big)<(K_1+K_2)(1+\log \K)<\infty .
\end{equation*}
We choose a natural number $n^*$ such that  
\begin{equation*}\label{perdida}
\sum_{n\geq n^*}\max \big(d_n^0,d_n^1\big)\leq \frac{\delta}{4},
\end{equation*}
where $\delta$ is the initial width of the strip for the  $\theta$ variable. The proof of Theorem \ref{teo} will be made using an algorithm divided in stages. The number  $n^*$ represents the starting stage for our algorithm. At each stage the functions will be analytic on a strip (for the $\theta$ variable) whose width equals 
\begin{eqnarray}
\delta^{n^*}&=&\frac{3\delta}{4}\\\label{bandita}
\delta^{n+1}&=&\delta^n-\max(d_n^0,d_n^1).
\end{eqnarray}   
The choice of  $n^*$ implies that for every  $n\geq n^*$ one has $\delta^n>\frac{\delta}{2}$. 
 \subsection{Description of the algorithm}\label{siete}
  The proof of Theorem  \ref{teo} will be obtained by showing that there exists an analytic conjugacy of the form 
  \begin{equation*}
  h(\theta,z)=\big(\theta,u_0(\theta)+u_1(\theta)z\big) 
  \end{equation*}
that puts one of the  $F_t$ into the adequate normal form,  with $u_1$ having zero degree. (As noted in Section \ref{FORMNORM}, this implies the desired result.) After conjugacy one obtains 
		\begin{displaymath}\begin{array}{lll}
 	h^{-1}\circ F\circ h&=&\Big (\theta+\alpha, \frac{1}{u_1(\tal)} \left  \{ \rho_{0,t}(\theta)+\lambda u_0(\theta)-u_0(\tal)\right \}\\
		&&+\quad \frac{1}{u_1(\tal)}\left \{ \rho_{1,t}(\theta)u_0(\theta)+\rho_t(\theta,u_0(\theta))\right \}\\
		&&+\quad z\frac{u_1(\theta)}{u_1(\tal)}\left \{ \lambda+ \rho_{1,t}(\theta)+\partial_z\rho_t(\theta,u_0(\theta))\right \}\\
&&+\quad \frac{1}{u_1(\tal)}\left \{ \rho_t(\theta,u_0(\theta)+u_1(\theta)z)-\rho_t(\theta,u_0(\theta))\right.\\
	&&\left. -\quad zu_1(\theta)\partial_z\rho_t(\theta,u_0(\theta))\right \}\Big ).
	\end{array}
	\end{displaymath}
	We will perform a process in order to eliminate the terms $\rho_0,\rho_1$ (this implies the desired normal form). This process is divided in many \emph{stages}. The goal of each stage is to decrease the sizes of $\rho_0$ and $\rho_1$, in such a way that at the limit these sizes are $0$. However, the rate of convergence is very slow and deeply depends on the arithmetical conditions for the pair  $(\alpha,\beta)$. 
  
  At each stage we will apply four  types of operations. The first two ones are KAM type (or Newton algorithm type) operations. That is, we solve the related equations but in a simplified way (we only keep terms which turns these equations into  cohomological equations). These operations  provide the functions  $u_0,u_1$  as solutions to cohomological equations, which give the corresponding coordinate change related to the stage. Obviously, this coordinate change does not conjugate  our transformation into the desired normal form, since we only solve simplified equations. We must deal with some \emph{rests},  whose sizes are controlled via the remaining two operations.  The third operation consists on  introducing a loss in the width of the analyticity strip for the $\theta$ variable. The fourth operation consists on restricting the parameter space. 
  
  We describe now the first two KAM type operations. These operations consist on solving  a truncate cohomological equation (cf. Section \ref{cohomo}): if we are in the stage $n$ we solve the equations with truncation up to the order  $2^n$.  
 \subsubsection{The four operations}     
 \paragraph{Case $u_1=1$ (\emph{solving the  $u_0$} equation). } After conjugacy by a map $h$ for which $u_1\equiv 1$ one obtains 
  \begin{displaymath}\begin{array}{lll}
 	h^{-1}\circ F\circ h&=&\Big (\theta+\alpha,  \rho_{0,t}(\theta)+\lambda u_0(\theta)-u_0(\tal) +  \rho_{1,t}(\theta)u_0(\theta)+\rho_t(\theta,u_0(\theta))\\
		&&+\quad z\left \{ \lambda+ \rho_{1,t}(\theta)+\partial_z\rho_t(\theta,u_0(\theta))\right \}\\
&&+\quad \left \{ \rho_t(\theta,u_0(\theta)+z)-\rho_t(\theta,u_0(\theta)) - z\partial_z\rho_t(\theta,u_0(\theta))\right \}\Big ) .
	\end{array}
	\end{displaymath}
Choosing $u_0$ so that 
\begin{equation}\label{ucero}
\rho_{0,t}(\theta)\Big|_{2^n}+\lambda u_0(\theta)-u_0(\tal)=0,
\end{equation}
the new dynamics is given by
\begin{eqnarray*}
\rho_{0,t}^{new}&=&\rho_{1,t}u_0+\rho_t(\cdot,u_0)+\rho_{0,t}\Big|_{\overline{2^n}}\\
\rho_{1,t}^{new}&=&\rho_{1,t}+\partial_z\rho_t(\cdot,u_0)\\
\rho_t^{new}(\cdot,z)&=&\rho_t(\cdot,u_0+z)-\rho_t(\cdot,u_0)-z\partial_z\rho_t(\cdot,u_0). \label{rhonew0}
\end{eqnarray*}
A crucial estimate will be 
\begin{equation}
\|u_0\|\leq C_22^n\Gamma_{\alpha,\beta}(2^n)\|\rho_{0,t}\|.
\end{equation}
Having in mind the sizes of $\rho_{0,t}$ and $\rho_{1,t}$, using this inequality we will see that the size of  $\rho^{new}_{0,t}-\rho_{0,t}\big|_{\overline{2^n}}$ decreases (``quadratic'' effect). The sizes of the remaining  involved functions remain almost unchanged.
\paragraph{Case $u_0=0$ (\emph{solving the  $u_1$} equation). } After conjugacy by a map $h$ for which $u_0\equiv 0$ one obtains
\begin{displaymath}\begin{array}{lll}
 	h^{-1}\circ F\circ h&=&\Big (\theta+\alpha, \frac{\rho_{0,t}(\theta)}{u_1(\tal)} + z\frac{u_1(\theta)}{u_1(\tal)}\left \{ \lambda+ \rho_{1,t}(\theta)\right \}\\
&&+\quad \frac{\rho_t(\theta,u_1(\theta)z)}{u_1(\tal)}\Big ).
		\end{array}
	\end{displaymath}
	It would be natural to choose $u_1$ so that 
\begin{equation}\label{lamdat}
\frac{u_1(\theta)}{u_1(\tal)}\left \{ \lambda+ \rho_{1,t}(\theta)\right \}=\tilde{\lambda}
\end{equation}
for an adequate  constant $\tilde{\lambda}$  (given by the zero mean condition in the cohomological equation). However, due to technical reasons, we will only deal with the corresponding 2n truncated cohomolical equation The resulting dynamics is then given by 
\begin{eqnarray*}
\rho_{0,t}^{new}&=&\frac{\rho_{0,t}(\cdot)}{u_1(\cdot+\alpha)}\\
\rho_{1,t}^{new}&=&\tilde{\lambda}-\lambda+Rest(\theta,t)\\\label{resta}
\rho_t^{new}&=& \frac{\rho_t(\cdot,u_1(\cdot)z)}{u_1(\cdot+\alpha)},\label{rhonew1}
\end{eqnarray*}
where the expression of  $Rest(\theta,t)$ will be made precise in the next Sections.  For  $u_1$ we will get an estimate of the form
\begin{equation}\label{estimau1d}
\|u_1-1\|\leq \K 2^n\ga(2^n)\|\partial_{\theta}\rho_{1,t}\|.
\end{equation}
The bound for the sizes of  $\rho_{0,t}^{new}$ and $\rho_t^{new}$ remains essentially unchanged.  The effect of solving an equation as the truncated version of (\ref{lamdat})  may be interpreted as  an attempt for turning  $\rho_1$ independent of  $\theta$. Indeed, this operation allows us to control the size of  $\partial_{\theta}\rho_{1,t}$.  However, to accomplish  this task we must control the size of  $Rest(\theta,z)$. 
 \paragraph{Loss in the width  of the strip. }The third operation consists on applying Cauchy estimates on the sizes of the functions. For this, we take the $\sup$ on a smaller strip. We call this operation a \emph{loss on the analyticity  strip}. By this way,  we control the resulting rest from  the two precedent operations. At the stage  $n$ we will lose  $\max (d_n^0,d_n^1)$. The new width  of the strip will be as in  (\ref{bandita}). This operation provides a control on the size of  $\rho_{0,t}^{new}$ and $\partial_{\theta}\rho_{1,t}^{new}$.
 \paragraph{Reduction of  the parameter space. }
Due to Claim  \ref{afirmac}, in order to control the size of  $\rho_{1,t}^{new}$ we need to control the size of  $\partial_{\theta}\rho_{1,t}^{new}$ and the complex number  $\int_{\T}\rho_{1,t}^{new}(\theta)d\theta$. The latter complex number may be changed with the parameter $t$ due to  the transversality condition on  the family $\{F_t\}_{t\in \Sigma}$. Thus, this operation consists on localizing a disc in the parameter space where this complex number is small enough. We pick  a disc centered at a simple zero of   $\int_{\T}\rho_{0,t}^{new}(\theta)d\theta$, and with a very small radius.
 \subsubsection{Description of a stage} 
 This Section only presents an informal description, and the complete details will be postponed until Section \ref{realisacion}. At the beginning of the  $n-$th stage we have the functions $\rho_{0,t}^n,\rho_{1,t}^n,\rho_t^n$, and the estimates 
 \begin{eqnarray}
\|\rho_{0,t}^{n}\|&\leq& w_n\\
\|\rho_{1,t}^{n}\|&\leq&K_{3}l_n
\end{eqnarray}
for a constant $K_{3}$. The sequences  $l_n,w_n$ are those defined in Section \ref{lemas}. We also have a parameter space  disc $D(t_n,p_n)$. The goal at the end of the stage is to achieve the corresponding necessary  estimates  to begin with the  $n+1$-th stage. Each stage splits into $4$ \emph{parts}:
\paragraph{First part. }  In this part we apply  many times the   ``solving the  $u_0$ equation'' operation. Hence,  this part is composed by  many \emph{steps}. We introduce a supplementary index indicating the corresponding step. We will obtain the functions  $\rho_{0,t}^{n,i+1}, \rho_{1,t}^{n,i+1},\rho_t^{n,i+1}$ as the result of the step $i$. At each step the sizes of the functions  $\rho_{1,t},\rho_t$ remain  essentially unchanged. Function  $\rho_{0,t}^{n,i}$ is the sum of a  \emph{quadratic} term whose size decrease in a half at each step, plus the series of rests having orders greater than  $2^n$ (Fourier series orders). The first part ends when the size of the quadratic term is smaller than  $\frac{w_{n+1}}{4}$. 
 \paragraph{Second part. } This part consists in applying one time the  ``solving the  $u_1$ equation'' operation. This enables  to relate the  control on the   size of  $\partial_{\theta}\rho_{1,t}$ to a control of a rest having order greater than  $2^n$.
 \paragraph{Third part. } This part consists in introducing a loss in the width of the strip. The first two  parts provide some rests arising essentially from truncate functions whose sizes are smaller than  $w_n, l_n$,  respectively. Inequality  (\ref{estimarestos}), the relations  (\ref{relad0}), (\ref{relad1}), and a loss on the strip width of an amount of  $\max (d_n^0,d_n^1)$, provide the estimates $\|\rho_{0,t}^{n+1}\|\leq w_{n+1}, \|\partial_{\theta}\rho_{1,t}^{n+1}\|\leq l_{n+1}$.
 \paragraph{Fourth part. } This part consists in  applying  the operation which reduces  the parameter space.  We get the bound $l_{n+1}$ for the complex number  $\int_{\T}\rho_{1,t}^{n+1}$.   Hence,  we  obtain all the necessary estimates to start with the $n+1$-th stage.

 \paragraph{ }
 Therefore, the proof of Theorem \ref{teo} depends on the possibility of showing that the hypotheses on the family  $\{F_t\}_{t\in \Sigma}$ allows, on the one hand, to obtain good estimates to start the first stage  (stage  $n^*$), and on the other hand, to complete our plan for any stage. 
 We will see in Section  \ref{convergencia} that the properties ensured by the Lemmas \ref{ele1} and  \ref{we1} for the sequences $l_n,w_n$ allow to show, for only one parameter  $t^*$ in $\Sigma$, the convergence of our process. In other words,  the successive composition of coordinates changes provided by each stage converges to an analytical coordinates change which  puts  the fibered holomorphic transformation  $F_{t^*}$  in an adequate  normal form.

\subsection{Initializing the algorithm}\label{depapropab}

At the beginning of our process we will deal with functions defined on the following domains:  the unitary disc  $\D=D(0,1)\subset \C$  for the  $z$ variable; the strip  $B_{\delta}$ of width  $\delta$ for the $\theta$ variable; and the disc  $D(0,K'\eps)\subset \C$ for the parameter  $t$, where the constant $K'$ will be explicit  very soon (see equation (\ref{radiot})).  
Before starting the iterative process we need to estimate the derivatives with respect to  $\theta$ of 
some functions. Allowing a loss of  $\delta/4$ in the width of the strip, classical Cauchy estimates give
\begin{eqnarray*}
\|\partial_{\theta}\rho_{1,t}\|&\leq& \frac{24\eps}{\delta}\\
\|\partial_{\theta}\partial_z\rho_t\|&\leq& \frac{24M}{\delta}=:N,
\end{eqnarray*} 
where the last equality provides the definition of the constant $N$. We may already exhibit the minimal value for the starting sizes, namely, 
\begin{equation}\label{epsilon}
\bar{\eps}=\frac{\delta \ l_{n^*}}{24}.
\end{equation}
These sizes only depend on  $\delta$ and $n^*$, which in their turn only depend on  $\hypos$. We will assume in this Section that  
 \begin{equation}
 \|\rho_0\|\leq w_{n^*}.
 \end{equation}
Notice that this hypothesis is stronger than those required  in the statement of Theorem  \ref{teo}. However, we will see that we can reduce the general case to this one by a preliminary process whose description we delay up to Section \ref{epsiloncuadrado}. We consider the functions 
 \begin{equation}
 \rho_{0,t}^{n^*,0}=\rho_{0,t}\quad ,\quad \rho_{1,t}^{n^*,0}=\rho_{1,t}\quad ,\quad \rho_t^{n^*,0}=\rho_t
 \end{equation}
defined on the following domain:  the strip of width  $\delta^{n^*}=3\delta/4$ for $\theta$, and the disc $D(0,R^{n^*,0})$ of radius  $R^{n^*,0}=1$ for   $z$. The parameter space is the disc  $D(t_{n^*},p_{n^*})$, centered at $t_{n^*}=0$. We will assume that $t_{n^*}$ is a simple zero of  $\int_{\T}\rho_{1,t}^{n^*,0}$.  (This hypothesis will be removed at  Section \ref{zerosimple}.)  The radius of the parameter space disc is  
 \begin{equation}
 p_{n^*}=100l_{n^*}=\frac{2400\bar{\eps}}{\delta}.
 \end{equation}
 Hence we pick   \begin{equation}\label{radiot}
 K'=\frac{2400}{\delta}.
 \end{equation} 
  Since $\K$ is a large number, we have  that $p_{n^*}<1$. Up so far, we have the following estimates to start with the stage $n^*$ :
	\begin{eqnarray}		
		\|\partial_{\theta}\rho_{1,t}^{n^*,0}\|&\leq&l_{n^*}\\
		\|\rho_{0,t}^{n^*,0}\|&\leq&w_{n^*}\\
	\|\partial_z^2\rho_t^{n^*,0}\|&\leq&M\\		osc\big(\partial_{\theta}\partial_z\rho_t^{n^*,0}\big)&\leq&2N\\
		osc\big(\partial_t\partial_z\rho_t^{n^*,0}\big)&\leq& 2T\\\label{esamisma}
				\Big\|\int_{\T}\partial_t\rho_{1,t}^{n^*,0}d\theta-\Delta_0\Big\|&\leq&\frac{L^{-1}}{1000}.
		 	\end{eqnarray}
			Notice that inequality  (\ref{esamisma}) and the existence of the complex number  $\Delta_0$ are also additional hypothesis which we will assume for a while. (These hypothesis will be removed in Section \ref{previous}.) Finally, we assume that  $L>|\Delta_0|>L^{-1}$.
			
\begin{rema}[On the starting stage and size of  $\eps$]\label{remarcaq} We will see in Section  \ref{remarcaepsilon} that a subtle modification in the first stage yields the last claim of Theorem  \ref{teo}, which  
concerns the sizes of  $u$ and $\eps$.
\end{rema}
 \subsection{Realization of the iterative scheme} \label{realisacion}
 To control the many functions arising during the realization of the iterative scheme, we introduce some real sequences. (The justification of some claims on these sequences made below will be deferred until they become transparent from the point of view of the iterative process.) At the begining of the stage $n$ we have the family  $\{F_t\}$ written in the form 
 \begin{equation*}
 F^{n,0}_t(\theta,z)=\big( \theta+\alpha,\rho_{0,t}^{n,0}(\theta)+z\rho_{1,t}^{n,0}(\theta)+\lambda z+\rho_t^{n,0}(\theta,z)\big ),
 \end{equation*}
where $\rho^{n,0}_t$ vanishes up to the order 2 at  $z=0$.  The definition domain for $F_t^{n,0}$ is the  disc  $D(0,R^{n,0})$ for the $z$ variable, where the sequence $\{R^{n,i}\}_{n\geq n^*,i\in\N\cup\{\infty\} }$ verifies
 \begin{equation}\label{Rdecrece}
 \frac{3}{8}<\dots<R^{n^*+1,1}<R^{n^*+1,0}<R^{n^*,\infty}<\dots <R^{n^*,2}<R^{n^*,1}<R^{n^*,0}=1,
 \end{equation}
and the strip  $B_{\delta^n}$ of width $\delta^n$ defined by (\ref{bandita}) for the $\theta$ variable. The parameter space is the disc $D(t_n,p_n)$ which is centered at   a simple zero $t_n$ of $\int_{\T}\rho_{1,t}^{n,0}(\theta)d\theta$, and whose  radius equals   
 \begin{equation}\label{PN}
 p_n=100l_n.
 \end{equation}
 In this Section we will deal only with the estimates which do not concern the parameter $t$, and we defer to  Section \ref{estimacionent} the estimates involving this parameter. At the beginning of the stage we have the estimates  
 \begin{eqnarray}\label{hiporo0}
 \|\rho_{0,t}^{n,0}\|&\leq& w_n\\
 \|\partial_z^2\rho_t^{n,0}\|&\leq& M_n\\
 osc\big(\partial_{\theta}\partial_z\rho_t^{n,0}\big)&\leq&N_{n,0}\\\label{dtetaro1}
 \|\partial_{\theta}\rho_{1,t}^{n,0}\|&\leq& l_n,
 \end{eqnarray}
 where the real sequences  $\{M_n\}_{n\geq n^*}$ and  $\{N_{n,i}\}_{n\geq n^*,i\in\N\cup\{\infty\}}$  are bounded from above by some constants  $K_M$ and $K_N$ respectively. 

Using Claim \ref{afirmac}, inequality  (\ref{dtetaro1}), and an estimate on the norm  $\big\|\int_{\T}\partial_{t}\rho_{1,t}^{n,0}-\Delta_0\big\|$ (see  inequality (\ref{hipofro1}) in Section \ref{estimacionent}), we obtain 
  \begin{equation}\label{hiporo1}
  \|\rho_{1,t}^{n,0}\|\leq  K_{3}l_n
  \end{equation}
   for some constant $K_3>0$. 
    \subsubsection{First part: the  $u_0$ equation}
  Suppose that at the first part of the stage  $n$ we have the functions $\rho_{0,t}^{n,i},\rho_{1,t}^{n,i},$ and $\rho_t^{n,i}$, where the last one is defined on the disc $D(0,R^{n,i})$ for the  $z$ variable. We will describe the $i$-th step. In this step we search  for a coordinate change of the form  $h_{n,i}^t(\theta,z)=(\theta,u_{0,t}^{n,i}(\theta)+z)$ by solving the truncate equation
  \begin{equation}\label{ucerotro}
  \lambda u^{n,i}_{0,t}(\theta)-u_{0,t}^{n,i}(\theta+\alpha)=-\rho_{0,t}^{n,i}|_{2^n}.
\end{equation}
By the results of Section \ref{cohomo}, the solution $u_{0,t}^{n,i}$ satisfies
\begin{eqnarray}\label{culo1}
\|u_{0,t}^{n,i}\|&\leq& C_22^n\Gamma_{\alpha,\beta}(2^n)\|\rho_{0,t}^{n,i}\|\\\label{culo2}
 \|\partial_{\theta}u_0^{n,i}\|&\leq& C_24^n\Gamma_{\alpha,\beta}(2^n)\|\rho_{0,t}^{n,i}\|.  \end{eqnarray}
Notice that the resulting new $\rho_0$ is
\begin{equation}\label{newrhocero000}
\rho_{0,t}^{n,i+1}=\rho_{0,t}^{n,i}|_{\overline{2^n}}+\rho_{1,t}^{n,i}u_{0,t}^{n,i}+\rho_t^{n,i}(\cdot,u_{0,t}^{n,i}),
\end{equation}
which due to the definition of $\rho_{0,t}^{n,i}|_{\overline{2^n}}$ satisfies   
\begin{equation}\label{tronco}
\rho_{0,t}^{n,i+1}|_{2^n}=\big\{\rho_{1,t}^{n,i}u_{0,t}^{n,i}+\rho_t^{n,i}(\cdot,u_{0,t}^{n,i})\big\}\big|_{2^n}.
\end{equation}
We also get a new  $\rho_1$, 
\begin{equation}\label{newro}
\rho_{1,t}^{n,i+1}=\rho_{1,t}^{n,i}+\partial_z\rho_t^{n,i}(\cdot,u_0^{n,i}),
\end{equation}
as well as a new $\rho$,
\begin{equation}\label{rho0}
\rho_t^{n,i+1}(\cdot,z)=\rho_t^{n,i}(\cdot,u_{0,t}^{n,i}+z)-\rho_t^{n,i}(\cdot,u_{0,t}^{n,i})-z\partial_z\rho_t^{n,i}(\cdot,u_{0,t}^{n,i}).
\end{equation}
The function $\_t^{n,i}(\theta,\cdot)$ is defined on  $D(0,r_t^{n,i+1}(\theta))$, with
\begin{equation}\label{erre}
r_t^{n,i+1}(\theta)= R^{n,i}-|u_{0,t}^{n,i}(\theta)|.
\end{equation}
The sequence $\{R^{n,i}\}$ verifies
\begin{equation}\label{erre0}
R^{n,i+1}\leq r_t^{n,i+1}(\theta)
\end{equation}
 for every $\theta$ in $B_{\delta^n}$. From now on the domain for $z$ is $D(0,R^{n,i+1})$.  
We recall that equation (\ref{ucerotro}) is solved without lossing width on the analyticity strip. In fact, we will introduce a loss of width of size $d_0^n$ only at the moment when we will need to estimate the rests  $\rho_{0,t}^{n,i}|_{\overline{2^n}}$. This will be made only once per stage, and at the third part of the each stage.  
\paragraph{Estimates at the first part of the stage.}
We define
\begin{eqnarray}
\eta_{n,0}^t&=&\rho_{0,t}^{n,0}\\
\eta_{n,i+1}^t&=&\rho_{1,t}^{n,i}u_{0,t}^{n,i}+\rho_t^{n,i}(\theta,u_{0,t}^{n,i}).
\end{eqnarray}
Since  (\ref{tronco}) yields 
\begin{equation}
\rho_{0,t}^{n,i+1}|_{2^n}=\eta_{n,i+1}^t|_{2^n},
\end{equation}
each time that we solve equation (\ref{ucerotro}) we actually obtain a better bound on  $\eta_{n,i+1}^t$.
This allows to improve (\ref{culo1}) and (\ref{culo2}) to
\begin{eqnarray}\label{culo10}
\|u_{0,t}^{n,i}\|&\leq& C_22^n\Gamma_{\alpha,\beta}(2^n)\|\eta_{n,i}^t\|_\\
\label{culo20}
\|\partial_{\theta}u_{0,t}^{n,i}\|&\leq& C_24^n\Gamma_{\alpha,\beta}(2^n)\|\eta_{n,i}^t\|.
\end{eqnarray}
We will show by induction on the index $i\geq 0$ of the corresponding step that, 
for some constant $K_4 \geq K_3$,   
\begin{eqnarray}\label{iM}
\|\partial_z^2\rho_t^{n,i}\|&\leq& M_n\\\label{iN}
osc \big(\partial_{\theta}\partial_z\rho_t^{n,i}\big)&\leq& N_{n,i}\\\label{iW}
\|\eta_{n,i}^t\|&\leq&\frac{w_n}{2^i}\\\label{iL}
\|\partial_{\theta}\rho_{1,t}^{n,i}\|&\leq& 2l_n\\\label{iRO}
\|\rho_{1,t}^{n,i}\|&\leq&K_4l_n.
\end{eqnarray}
During the proof we will also provide an inductive definition for the sequences $\{R^{n,i}\}_{i\geq 0}$ and $\{N_{n,i}\}_{i\geq 0}$. 
Notice that inequalities  (\ref{culo10}), (\ref{culo20}) and (\ref{iW}) give
\begin{eqnarray}\label{culo3}
\|u_{0,t}^{n,i}\|&\leq& \frac{C_2l_n}{\K2^n2^i}\\\label{culo4}
\|\partial_{\theta}u_{0,t}^{n,i}\|&\leq& \frac{C_2l_n}{\K2^i}.
\end{eqnarray}
The case $i=0$ for the induction is furnished by the hypothesis (\ref{hiporo0})$, \dots, $(\ref{hiporo1}). We now assume that the desired relations hold for every  $j\in \{0,\dots,i\}$ and that we have at our disposal the constants $R^{n,j}$ and $N_{n,j}$  for every  $j\in \{0,\dots,i\}$. We let 
\begin{eqnarray}
R^{n,i+1}&=&R^{n,i}-C_22^n\Gamma_{\alpha,\beta}(2^n)w_n\Big(\frac{1}{2^i}\Big)\\
R^{n,\infty}&=&R^{n,0}-2C_22^n\Gamma_{\alpha,\beta}(2^n)w_n.
\end{eqnarray}
If $\K$ is larger than  $C_2$ then inequality (\ref{erre0}) follows from (\ref{erre}) and (\ref{culo3}). 
Inequality (\ref{iRO}) follows from (\ref{newro}), (\ref{hiporo1}), (\ref{iM}) and (\ref{culo3}), since
\begin{eqnarray*}
\|\rho_{1,t}^{n,i+1}\|&\leq& \bigg\|\rho_{1,t}^{n,0}+\sum_{j=0}^i\partial_z\rho_t^{n,j}(\theta,u_{0,t}^{n,j})\bigg\|\\
&\leq&K_{3}l_n+M_n\sum_{j=0}^i\|u_{0,t}^{n,j}\|\\
&\leq&K_{3}l_n+K_MC_2\frac{l_n}{2^n}\\
&\leq&K_4l_n.
\end{eqnarray*}
Notice that the constant $K_4$ only depends on $C_2,K_3$ and $K_M$. From relation (\ref{rho0}) we obtain 
\begin{equation*}
\|\partial_z^2\rho_t^{n,i+1}\|=\|\partial_z^2\rho_t^{n,i}(\cdot,z+u_0^{n,i})\|\leq M_n,
\end{equation*}
which implies   (\ref{iM}). Differentiating    (\ref{rho0}) we get 
\begin{eqnarray*}
\partial_{\theta}\partial_z\rho_t^{n,i+1}&=&\partial_{\theta}\partial_z\rho_t^{n,i}(\cdot,z+u_{0,t}^{n,i})+\partial_z^2\rho_t^{n,i}(\cdot,z+u_{0,t}^{n,i})\partial_{\theta}u_{0,t}^{n,i}\\
&&-\partial_{\theta}\partial_z\rho_t^{n,i}(\cdot,u_{0,t}^{n,i})-\partial_z^2\rho_t^{n,i}(\cdot,u_{0,t}^{n,i})\partial_{\theta}u_{0,t}^{n,i}.\label{newrhotheta0}
\end{eqnarray*}
By estimating the corresponding terms we can bound the oscillation
\begin{eqnarray*}
osc \big(\partial_{\theta}\partial_z\rho_t^{n,i+1}\big)&\leq& N_{n,i}+2K_MK_24^n\Gamma_{\alpha,\beta}(2^n)\frac{w_n}{2^i}.
\end{eqnarray*}
The above  inequality allows us to give the definition of the sequence  $\{N_{n,i}\}_{i\geq 0}$: letting  $N_{n^*,0}=2N$ and assuming that $N_{n,i}$ is already defined, we put
\begin{equation}
N_{n,i+1}=N_{n,i}+2K_MC_24^n\Gamma_{\alpha,\beta}(2^n)\frac{w_n}{2^i}.\\
\end{equation}
We also define $N_{n,\infty}$ as the supremum of the $N_{n,i}$'s, that is,
\begin{equation}\label{ninfinito}
N_{n,\infty}=N_{n,0}+4K_MC_24^n\Gamma_{\alpha,\beta}(2^n)w_n.
\end{equation}
It is not hard to check that  
\begin{equation*}
osc \big(\partial_{\theta}\partial_z\rho_t^{n,i+1}\big)\leq N_{n,i+1}.
\end{equation*}
Relation (\ref{newro}) allows us to (compute and) estimate  $\partial_{\theta}\rho_1^{n,i+1}$ by 
\begin{eqnarray*}
\|\partial_{\theta}\rho_{1,t}^{n,i+1}\|&=&\bigg \|\partial_{\theta}\rho_{1,t}^{n,0}+\sum_{j=0}^i\partial_{\theta}\Big(\partial_z\rho_t^{n,j}(\theta,u_{0,t}^{n,j})\Big)\bigg \|\\
&\leq&l_n+\bigg\|\sum_{j=0}^i\partial_{\theta}\partial_z\rho_t^{n,j}(\theta,u_{0,t}^{n,j})+\partial^2_z\rho_t^{n,j}(\theta,u_{0,t}^{n,j})\partial_{\theta}u_{0,t}^{n,j}\bigg \|\\
&\leq&l_n+\sum_{j=0}^iN_{n,j}\|u_{0,t}^{n,j}\|+M_n\|\partial_{\theta}u_{0,t}^{n,j}\|\\
&\leq&l_n+\frac{K_5}{\K}l_n,
\end{eqnarray*}
where $K_5$ is a constant which only depends on $C_2,K_M$ and $K_N$. If we pick $\K$ larger than $K_5$ we get
\begin{equation*}
\|\partial_{\theta}\rho_{1,t}^{n,i+1}\|\leq 2l_n.
\end{equation*}
Finally, we estimate $\eta_{n,i+1}^t$ using point \emph{3.} of Lemma  \ref{ele1} and point \emph{1.} of Lemma \ref{we1},
\begin{eqnarray*}\label{etaa}
\|\eta_{n,i+1}^t\|&\leq& K_4l_nC_22^n\Gamma_{\alpha,\beta}(2^n)\frac{w_n}{2^i}+M_n\bigg (C_22^n\Gamma_{\alpha,\beta}(2^n)\frac{w_n}{2^i}\bigg)^2\\
&\leq& \frac{1}{16}\frac{w_n}{2^i}+\frac{1}{2^i16}\frac{w_n}{2^i}\\
&\leq&\frac{w_n}{2^{i+1}}.
\end{eqnarray*}
We used here that $\K>K_4C_2$ and $\K>K_MC_2^2$.
\paragraph{End of the first part. } We apply the operation ``solving the  $u_0$ equation'' until the size of  $\|\eta_{n,i}^t\|$ becomes smaller than  $\frac{w_{n+1}}{4}$; more precisely, when  $\frac{w_n}{2^i}$ is smaller than $\frac{w_{n+1}}{4}$. At this time $i_n$ we stop and we pass to the second part of the stage.
\subsubsection{Second part: the $u_1$ equation}
This part consists in solving the equation for $u_1$ only once. More precisely, we search for a coordinates  change of the form $h_n^t(\theta,z)=(\theta,\tilde u^n_{1,t}(\theta)z)$ such that 
\begin{equation}\label{star}
\frac{\tilde u_{1,t}^n(\theta)}{\tilde u_{1,t}^n(\theta+\alpha)}\big \{ \rho_{1,t}^{n,i_n}(\theta)+\lambda\big \}=\lambda_n^t.
\end{equation}
Notice that if we put  $e^{\tilde v^n_t(\theta)}=\tilde u_{1,t}^n(\theta)$, equation (\ref{star}) reduces to solving  
 \begin{equation}
 \tilde v^n_t(\theta)-\tilde v^n_t(\theta+\alpha)=\log \lambda_n^t-\log(\rho_{1,t}^{n,i_n}(\theta)+\lambda).
\end{equation}
However, due to technical reasons we will just deal with the following truncate equation
\begin{equation}\label{uuntro}
v^n_t(\theta)-v^n_t(\theta+\alpha)=\log \lambda_n^t-\big\{\log\big(\rho_{1,t}^{n,i_n}(\theta)+\lambda \big)\big\}\big|_{2^n}.
\end{equation}
The complex number  $\lambda_n^t$ is  the unique one for which a continuous solution to this equation exists (see Section \ref{cohomo}). More precisely, it corresponds to the unique complex number for which  the mean of the r.h.s. vanishes,  that is, 
\begin{equation}\label{valorlamda}
\log \lambda_n^t=\int_{\T}\big\{\log\big (\rho_{1,t}^{n,i_n}(\theta)+\lambda\big )\big \}\big|_{2^n}d\theta=\log \lambda +\int_{\T}\log \Big(\frac{\rho_{1,t}^{n,i_n}}{\lambda}+1\Big)d\theta.
\end{equation}
We have 
\[
\Big |\int_{\T}\log \Big(\frac{\rho_1^{n,i_n}}{\lambda}+1\Big)d\theta \Big|\leq C_3\|\rho_1^{n,i_n}\|
\]
for a constant   $C_3$. To estimate the solution of  (\ref{uuntro}) we need an estimate on the size of the r.h.s. 
For this, we rewrite this expression as
\begin{equation*}
\log \lambda_n^t -\log \big ( \rho_{1,t}^{n,i_n}(\theta)+\lambda\big )=\int_{\T}\log \Big(\frac{\rho_{1,t}^{n,i_n}}{\lambda}+1\Big)d\theta - \log \Big(\frac{\rho_{1,t}^{n,i_n}}{\lambda}+1\Big)
\end{equation*}
and we obtain
\begin{equation*}
\|\log \lambda_n^t -\log \big ( \rho_{1,t}^{n,i_n}(\theta)+\lambda\big )\|\leq C_4\|\rho_{1,t}^{n,i_n}\|
\end{equation*}
for a constant  $C_4$. We have
\begin{eqnarray}\label{vvv}
 \|v^n_t\|&\leq& C_2C_42^n\Gamma_{\alpha}(2^n)\|\rho_{1,t}^{n,i_n}\|\\\label{dv}
 \|\partial_{\theta}v^n_t\|&\leq& C_2C_44^n\Gamma_{\alpha}(2^n)\|\rho_{1,t}^{n,i_n}\|
 \end{eqnarray}
 since  we choose the zero mean solution. Letting $u_{1,t}^n=\exp(v_{1,t}^n)$ and conjugating by  $h_n^t=(\theta,u_{1,t}^n(\theta)z)$ we get the new functions 
 \begin{eqnarray}\label{rhoceronew}
 \rho_{0,t}^{n+1,0}&=&\frac{\rho_{0,t}^{n,i_n}(\cdot)}{u_{1,t}^n(\cdot+\alpha)}\\ \label{rhoesnew}
 \rho_t^{n+1,0}&=&\frac{1}{u_{1,t}^n(\cdot+\alpha)}\rho_t^{n,i_n}(\cdot,u_{1,t}^nz).
 \end{eqnarray}
 The last one is defined on a disc  $D(0,r_t^{n+1,0}(\theta))$ for the  $z$ variable, with 
 \begin{equation}\label{puro}
 r_t^{n+1,0}(\theta)=\frac{R^{n,i_n}}{|u_{1,t}^n(\theta)|}.
 \end{equation}
 The sequence $\{R^{n,i}\}$ verifies 
 \begin{equation}\label{erre1nuevo}
 R^{n+1,0}\leq r_t^{n+1,0}(\theta)
 \end{equation}
 for every $\theta$ in the strip $B_{\delta^n}$.  From now on,  the domain for the $z$ variable is the disc  $D(0,R^{n+1,0})$. 
 The untouched rest in equation  (\ref{uuntro}) gives raise to a new function $\rho_1$ defined by
 \begin{equation}\label{rhounonew}
  \rho_{1,t}^{n+1,0}(\theta)=\lambda_n^t\Big ( e^{\big\{ \log(\rho_{1,t}^{n,i_n}(\theta)+\lambda)\big\}\big|_{\overline{2^n}}}\Big ) -\lambda .\\ 
 \end{equation}
 \paragraph{Estimates at the second part of the stage.} By  (\ref{vvv}) and (\ref{iRO}) we have
 \begin{equation}\label{Kocho}
 \|u_{1,t}^n-1\|\leq C_3C_2C_42^n\Gamma_{\alpha}(2^n)K_4l_n.
 \end{equation}
Therefore, there exists a constant   $K_6$   such that
 \begin{eqnarray}\label{uuno}
 \|u_{1,t}^n\|&\leq& e^{K_62^n\Gamma_{\alpha}(2^n)l_n}\\ \label{duuno}
 \|(u_{1,t}^n)^{-1}\|&\leq& e^{K_62^n\Gamma_{\alpha}(2^n)l_n}.
  \end{eqnarray}
  Using the equality   $\partial_{\theta}u_1^n=\partial_{\theta}v^nu_1^n$ and (\ref{dv})  we obtain
  \begin{eqnarray}
 \|\partial_{\theta}u_{1,t}^n\|&\leq& e^{K_62^n\Gamma_{\alpha}(2^n)l_n}K_64^n\Gamma_{\alpha}(2^n)l_n\\\label{detetau}
 &\leq& 3K_64^n\ga(2^n)l_n
   \end{eqnarray}
   provided that  $\K$ is larger than $K_6$.   
      We remark that the estimate  (\ref{Kocho}) says that the topological degree of  $u_1$ is zero. We may now proceed to conclude the definition of the sequence $\{R^{n,i}\}$ and to prove the lower bound $3/8$ for it.
   \paragraph{The sequence $\{R^{n,i}\}$.} Assume that  $R^{n,\infty}$ is already defined. We put  
   \begin{equation}
   R^{n+1,0}=\frac{R^{n,\infty}}{e^{\K2^n\Gamma_{\alpha}(2^n)l_n}}.
   \end{equation}
   With this definition,  inequality  (\ref{erre1nuevo}) follows from   (\ref{puro}) and (\ref{duuno}). The sequence  $R^{n,i}$ also satisfies the monotonicity  properties announced in   (\ref{Rdecrece}). Finally, to check 
   that  $R^{n,0}$ (and therefore each $R^{n,i}$) is bounded from below by $\frac{3}{8}$, we first note that  
   \begin{equation*}
R^{n+1,0}=1-\sum_{j=0}^{n} \left(R^{j,0}-R^{j+1,0}\right ).
\end{equation*}
Each of the differences above can be estimated by 
\begin{eqnarray*}
R^{j,0}-R^{j+1,0}&=&R^{j,0}( 1-e^{-\K2^j\Gamma_{\alpha}(2^j)l_j}) + \frac{2C_22^j\Gamma_{\alpha,\beta}(2^j)w_j}{e^{\K2^j\Gamma_{\alpha}(2^j)l_j}}\\
&\leq& \K2^j\Gamma_{\alpha}(2^j)l_j+2C_22^j\Gamma_{\alpha,\beta}(2^j)w_j.
\end{eqnarray*}
Thus, from point  \emph{2.} of Lemma  \ref{ele1} and point  \emph{2.} of Lemma  \ref{we1},
\begin{equation*}
\sum_{0}^n\left (R^{j,0}-R^{j+1,0}\right)\leq \sum_n\K2^n\Gamma_{\alpha}(2^n)l_n + 2\sum_n C_22^n\Gamma_{\alpha,\beta}(2^n)w_n<\frac{5}{8}\quad_{\blacksquare}
\end{equation*}
\newline

Relation  (\ref{rhoesnew}) gives 
 \begin{equation*}
 \partial_z^2\rho_t^{n+1,0}(\theta,z)=\frac{u_{1,t}^n(\theta)^2}{u_{1,t}(\theta+\alpha)}\partial_z^2\rho_t^{n,i_n}(\theta,u_{1,t}^n(\theta)z).
 \end{equation*}
 Therefore, according to the definition below of the sequence $\{M_n\}$ we get the  desired estimate for $\partial_z^2\rho_t$, namely, 
 \begin{eqnarray}\nonumber
 \|\partial_z^2\rho_t^{n+1,0}\|&\leq& M_ne^{3K_62^n\ga(2^n)l_n}\\\label{d2ronuevo}
 &=& M_{n+1}.
 \end{eqnarray}
 \paragraph{The sequence  $\{M_n\}$.} We first let $M_{n^*}=M$, and assuming that $M_n$ is defined, we let
 \begin{equation}
 M_{n+1}=M_ne^{3\K 2^n\ga(2^n)l_n}.
 \end{equation}
 According to this definition, one easily checks the existence of an upper bound $K_M$ for the sequence $M_n$ $\quad_{\blacksquare}$
 \newline
 
 In order to compute the oscillation  $\partial_{\theta}\partial_z\rho_t^{n+1,0}$ we first note that
 \begin{eqnarray*}
 \partial_{\theta}\partial_z\rho_t^{n+1,0}(\theta,z)&=&\frac{\partial_{\theta}u^n_{1,t}(\theta)u^n_{1,t}(\theta+\alpha)-\partial_{\theta}u_{1,t}^n(\theta+\alpha)u_{1,t}^n(\theta)}{u_{1,t}^n(\theta+\alpha)^2}\partial_z\rho_t^{n,i_n}(\theta,u_{1,t}^n(\theta)z)\\
 &&+\frac{u_{1,t}^n(\theta)}{u_{1,t}^n(\theta+\alpha)}\Big\{\partial_{\theta}\partial_z\rho_t^{n,i_n}(\theta,u_{1,t}^n(\theta)z)\\
 &&+\partial_z^2\rho_t^{n,i_n}(\theta,u_{1,t}^n(\theta)z)z\partial_{\theta}u_{1,t}^n(\theta)\Big \},\label{newrhotheta1}
 \end{eqnarray*}
 which implies
  \begin{eqnarray}\nonumber
 osc \big(\partial_{\theta}\partial_z\rho_t^{n+1,0}\big)&\leq&4\|\partial_{\theta}u_{1,t}^n\|\|u_{1,t}^n\|\|(u_{1,t}^n)^{-1}\|^2K_M+N_{n,i_n}\|u_{1,t}^n\|\|(u_{1,t}^n)^{-1}\|\\\label{dtetarotnew}
 &&+2\|(u_{1,t}^n)^{-1}\|\|u_{1,t}^n\|K_M\|\partial_{\theta}u_{1,t}^n\|.
 \end{eqnarray}
At this point we can conclude the definition of the sequence $\{N_{n,i}\}$ and show the existence of an upper bound $K_N$ for it.
 \paragraph{The sequence  $\{N_{n,i}\}$.}  Assuming that  $N_{n,\infty}$ is already defined, we put 
 \begin{equation}
 N_{n+1,0}=e^{2K_62^n\ga(2^n)l_n}N_{n,\infty}+144K_MK_64^n\ga(2^n)l_n.
 \end{equation}
This definition together with  inequality  (\ref{dtetarotnew}) and other estimates in this Section give us 
\begin{equation}
osc \big(\partial_{\theta}\partial_z\rho_t^{n+1,0}\big)\leq N_{n+1,0}.
\end{equation}
 By replacing the value of  $N_{n,\infty}$ (see (\ref{ninfinito})) in the definition of  $N_{n+1,0}$ we obtain
 \begin{eqnarray*}
N_{n+1,0}&\leq& e^{2K_62^n\ga(2^n)l_n}N_{n,0}+8K_MC_24^n\gab(2^n)w_n+144K_MK_64^n\ga(2^n)l_n\\
&\leq&e^{2K_62^n\ga(2^n)l_n}N_{n,0}+K_74^n\ga(2^n)l_n\\\label{K7}
&\leq&e^{2K_6\sum2^n\ga(2^n)l_n}\Big(2N+K_7\sum4^n\ga(2^n)l_n \Big).
 \end{eqnarray*}
 These last inequalities hold for a well defined constant  $K_7$. We see that  (\ref{K7}) and the properties of the sequence  $\{l_n\}$ imply an upper bound $K_N$ for the  sub-sequence $N_{n,0}$.  Finally, due to the following monotonicity relations, 
\begin{equation*}\label{eneacotado}
2N=N_{n^*,0}<N_{n^*,1}<\dots <N_{n^*,\infty}<N_{n^*+1,0}<\dots
\end{equation*}
this implies the existence of an upper bound for the whole sequence $\{N_{n,i}\}$ $\quad_{\blacksquare}$
\newline

We conclude this second part of the stage by writing out the functions whose expressions involve some truncate functions and for which we need to allow a loss in the width of the strip  in order to get good estimates. First of all, relation  (\ref{rhounonew}) implies 
 \begin{equation*}
 \partial_{\theta}\rho_{1,t}^{n+1,0}=\lambda_n^t\Big ( e^{\big\{ \log(\rho_{1,t}^{n,i_n}+\lambda)\big\}\big|_{\overline{2^n}}}\Big )\partial_{\theta}\Big(\Big\{\log (\rho_{1,t}^{n,i_n}+\lambda)\Big \}\Big |_{\overline{2^n}} \Big).
 \end{equation*}
 Since the truncation operator commutes with the derivation with respect to $\theta$, we have 
 \begin{equation}\label{resto1}
 \partial_{\theta}\rho_{1,t}^{n+1,0}=\lambda_n^t\Big ( e^{\big\{ \log(\rho_{1,t}^{n,i_n}+\lambda)\big\}\big|_{\overline{2^n}}}\Big )\Big\{ \frac{1}{\rho_{1,t}^{n,i_n}+\lambda}\big ( \partial_{\theta}\rho_{1,t}^{n,i_n}\big )\Big \}\Big |_{\overline{2^n}}.
 \end{equation}
 Finally, the expression for  $\rho_0$ is 
\begin{equation}\label{resto0}
\rho_{0,n}^{n+1,0}=\frac{\rho_{0,t}^{n,i_{n}}}{u_{1,t}^n(\cdot+\alpha)}=\frac{1}{u_{1,t}^n(\cdot+\alpha)}\bigg(\eta_{n,i_n}^t+\sum_{j=0}^{i_n-1}\eta_{n,j}^t|_{\overline{2^n}}\bigg).
\end{equation}
\subsubsection{Third part: loss on the width of the strip}
In order to estimate the size of the truncate functions appearing in (\ref{resto1}) and (\ref{resto0}), we introduce a loss in the width of the strip of size  $\max(d_n^0,d_n^1)$. The new domain for the $\theta$ variable is $\delta^{n+1}$.  We can estimate  $\rho_{0,t}^{n+1,0}$ by
\begin{eqnarray*}
\|\rho_{0,t}^{n+1,0}\|&\leq& 2\bigg (\frac{w_{n+1}}{4}+\sum_{j=0}^{i_n-1}C_1\frac{w_ne^{-2\pi 2^n d_n^0}}{2^jd_n^0}\bigg)\\
&\leq& \frac{w_{n+1}}{2}+4C_1\frac{w_ne^{-2\pi 2^n d_n^0}}{d_n^0}\\\label{ultimawea}
&\leq& w_{n+1}
\end{eqnarray*}
provided that   $\K$ is larger than $2C_1$. Note that here we used the estimate $\|(u_{1,t}^n)^{-1}\|<2$. This estimate is not sharp but it is enough for our purposes.  To estimate the expression in (\ref{resto1}) we first notice that, since  $\|\rho_{1,t}^{n,i_n}\|$ is uniformly small, there exists  a constant  $C_5$ such that  
  \begin{eqnarray}\nonumber
 \|\partial_{\theta}\rho_{1,t}^{n+1,0}\|&\leq& \frac{C_5\|\partial_{\theta}\rho_{1,t}^{n,i_n}\|e^{-2\pi 2^n d_n^1}}{d_n^1}\\\nonumber
 &\leq&\frac{2C_5l_ne^{-2\pi 2^nd_n^1}}{d_n^1}\\\label{dtetanuevoro1}
 &\leq&l_{n+1}.
  \end{eqnarray} 
This last inequality holds thanks to the definition of the sequence  $\{l_n\}_{n\geq 0}$ (see Lemma \ref{ele1}) and the fact that  $\K$ is greater than  $2C_5$.  In this way, we get the family of inequalities (\ref{hiporo0})$\dots$(\ref{dtetaro1}) for the index  $n+1$.  This finishes the third part. 
  \subsection{Estimates with respect to the parameter}\label{estimacionent}
  Before dealing  with the fourth and last part of the stage, we treat again the first three parts in order to get some estimates  on the parameter $t$. These estimates will be essential in the fourth part, where we will reduce the parameter space in an important way. However, in this Section we will also reduce  the parameter space in order to get some useful Cauchy estimates. At the beginning of the stage $n$ we have the following bounds:
 \begin{eqnarray}\label{dtdzrho}
   osc\big(\partial_t\partial_z\rho_t^{n,0}\big)&\leq&T_{n,0}\\\label{hipofro1}
 \Big\|\int_{\T}\partial_t\rho_{1,t}^{n,0}d\theta-\Delta_0\Big\|&\leq& s_n .
 \end{eqnarray} 
 The sequence of real numbers  $\{T_{n,i}\}_{n\geq n^*,\in \N\cup\{\infty\}}$ verifies 
  \begin{equation}\label{mayT}
 T_{n+1,0}\leq K_Tn2^n\ga(2^n) 
\end{equation}
for a constant $K_T$. The sequence of real numbers $\{s_n\}_{n\geq n^*}$ verifies 
\begin{equation}\label{vales}
\frac{L^{-1}}{1000}=s_{n^*}<s_{n^*+1}<s_{n^*+2}<\dots<\frac{L^{-1}}{100}.
\end{equation}
 \paragraph{Estimates at the first part of the stage.} We will show inductively on $i\geq 0$ that 
 \begin{equation}
osc \big(\partial_t\partial_z\rho_t^{n,i}\big)\leq T_{n,i} .
 \end{equation}
 Suppose that this estimate holds for  $0\leq j< i+1$. We start by estimating   $\partial_tu_{0,t}^{n,i}$. For this we use a Cauchy estimate, that is,  we lose a small quantity in the radius of the disc  $D(t_n,p_n)$ (the parameter space) in order to get an estimate for  $\partial_tu_{0,t}^{n,i}$. More precisely, we construct a sequence of radius $\{p_{n^-_i}\}_{i\in \N \cup \infty}$ defined by 
 \begin{eqnarray*}
 p_{n^-_{-1}}&=&p_n\\
 p_{n^-_{i}}&=&p_{n^-_{i-1}}-l_n\Big(\frac{2}{3}\Big)^i\\
 p_{n^-_{\infty}}&=&p_n-3l_n.
 \end{eqnarray*}
 This says that at each step we lose a very small quantity of radius of size $l_n\Big(\frac{2}{3}\Big)^i$ (we use this loss in the estimate  (\ref{nosequewa}) that follows). We note that at the limit, the total loss has size $3l_n$.  From  (\ref{culo3}) the Cauchy estimate gives 
 \begin{equation}\label{nosequewa}
 \|\partial_t u_{0,t}^{n,i}\|\leq \frac{6 C_2}{\K2^n}\Big(\frac{3}{4}\Big)^i.
 \end{equation}
 We can now estimate the oscillation of  $\partial_t\partial_z\rho$ by
 \begin{eqnarray*}
 \partial_t\partial_z\rho_t^{n,i+1}&=&\partial_t\partial_z\rho_t^{n,i}(\theta,u_{0,t}^{n,i}+z)+\partial_z^2\rho_t^{n,i}(\theta,u_{0,t}^{n,i}+z)\partial_tu_{0,t}^{n,i}\\\label{rhonewt0}
 &&-\partial_t\partial_z\rho_t^{n,i}(\theta,u_{0,t}^{n,i})-\partial_z^2\rho_t^{n,i}(\theta,u_{0,t}^{n,i})\partial_tu_0^{n,i}\\
 osc \big(\partial_t\partial_z\rho_t^{n,i+1}\big)&\leq&T_{n,i}+2M_n\|\partial_tu_{0,t}^{n,i}\|.
 \end{eqnarray*}
 This enables  us to  define the sequence $\{T_{n,i}\}_{i\geq 0}$ by letting  $T_{n^*,0}=2T$ and defining recursively 
 \begin{equation}
 T_{n,i+1}=T_{n,i} +\frac{12K_M C_2}{\K2^n}\Big(\frac{3}{4}\Big)^i.
 \end{equation}
 We define  $T_{n,\infty}$ by 
 \begin{equation}\label{Tinfinito}
 T_{n,\infty}=T_{n,0} +\frac{48K_M C_2}{\K2^n}.
 \end{equation} 
Notice that $T_{n,0}$ is an upper bound for the $T_{n,i}$, and  the  following inequality holds
\begin{equation*}\label{newdtrot}
osc \big(\partial_t\partial_z\rho_t^{n,i+1}\big)\leq T_{n,i+1}.
\end{equation*} 
  \paragraph{Estimates at the second part of the stage.} First,  we estimate the size of  $\partial_t u_1^n=u_1^n\partial_tv_n$.  By taking derivatives with respect to $t$ in both sides of (\ref{uuntro}) we obtain
  \begin{equation}
  \partial_t v^n_t(\theta)-\partial_t v^n_t(\theta+\alpha)=-\Big\{ \frac{\partial_t\rho_{1,t}^{n,i_n}}{\rho_{1,t}^{n,i_n}+\lambda}\Big \}\Big |_1^{2^n}.
  \end{equation}
Recall that the notation  $|_1^{2^n}$ means that we only keep the Fourier series terms  with orders  between $1$ and $2^n$. The best estimate that we can get for the right hand side function  (in his original untruncated form) is of the order of a constant. More  precisely, by introducing a loss in the radius of the disc $D(t_n,p_{n^-_{\infty}})$ of size  $l_n$, we get a Cauchy estimate 
  \begin{equation}\label{silvio}
  \Big\|\frac{\partial_t\rho_{1,t}^{n,i_n}}{\rho_{1,t}^{n,i_n}+\lambda}\Big \|\leq 12K_4 .
  \end{equation}
  The new parameter space disc is  $D(t_n,p_{n^-})$, where
  \begin{equation*}
  p_{n^-}=p_{n^-_{\infty}}-l_n=p_n-4l_n.
\end{equation*}
    Then we have
 \begin{equation}\label{deteuuno}
 \|\partial_tu_{1,t}^n\|\leq C_224K_4 2^n\pb(2^n).
 \end{equation}
 We may also estimate the oscillation of $\partial_t\partial_z\rho_t^{n+1,0}$ by
\begin{eqnarray}\nonumber
 \partial_{t}\partial_z\rho_t^{n+1,0}(\theta,z)&=&\frac{\partial_{t}u^n_{1,t}(\theta)u^n_{1,t}(\theta+\alpha)-\partial_{t}u_{1,t}^n(\theta+\alpha)u_{1,t}^n(\theta)}{u_{1,t}^n(\theta+\alpha)^2}\partial_z\rho_t^{n,i_n}(\theta,u_{1,t}^n(\theta)z)\\\nonumber
 &&+\frac{u_{1,t}^n(\theta)}{u_{1,t}^n(\theta+\alpha)}\Big\{\partial_t\partial_z\rho_t^{n,i_n}(\theta,u_{1,t}^n(\theta)z)\\\nonumber
  &&+\partial_z^2\rho_t^{n,i_n}(\theta,u_{1,t}^n(\theta)z)z\partial_tu_{1,t}^n(\theta)\Big \}\\\nonumber
 osc \big(\partial_t\partial_z\rho_t^{n+1,0}\big)&\leq& 4\|\partial_tu_{1,t}^n\|\|u_{1,t}^n\|\|(u_{1,t}^n)^{-1}\|^2M_n+\|u_{1,t}^n\|\|(u_{1,t}^n)^{-1}\|T_{n,i_n}\\\label{rhonewt1}
 &&+2\|u_{1,t}^n\|\|(u_{1,t}^n)^{-1}\|M_n\|\partial_tu_{1,t}^n\|.
  \end{eqnarray}
  We   introduce now  the definition of the sequence $\{T_{n,0}\}_{n\geq n^*}$ and show the claimed upper bound  (\ref{mayT}).
 \paragraph{The sequence  $\{T_{n,i}\}$.} Suppose  $T_{n,\infty}$ is defined. We put  
 \begin{equation}
T_{n+1,0}=\tilde{K}_T2^n\ga(2^n)+2T_{n,\infty},
 \end{equation} 
where the constant $\tilde{K}_T$ is given by $\tilde{K}_T=32K_MC_224K_4$. By looking at the inequality  (\ref{rhonewt1}) we obtain 
  \begin{equation*}
 osc \big(\partial_t\partial_z\rho_t^{n+1,0}\big)\leq T_{n+1,0}.
 \end{equation*}
 If we replace the value of $T_{n,\infty}$ (see (\ref{Tinfinito})) in the definition of  $T_{n+1,0}$ we get
 \begin{equation*}
 T_{n+1,0}\leq\tilde{K}_T2^n\ga(2^n)+2T_{n,0}+\frac{\tilde{K}_T}{2^n}.
 \end{equation*}
 This gives the desired estimate (in fact, we get a sharper estimate) $\quad_{\blacksquare}$
\paragraph{Estimates at the third part of the stage.}   We can estimate the distance between $\int_{\T}\partial_t\rho_{1,t}$ and $\Delta_0$ and show that this distance is small. In Section \ref{coup} this will allow us  to show that $\int_{T}\partial_t\rho_{1,t}(\theta)d\theta$ grows almost like a linear map. 
   In order to simplify the notation, we introduce the constants $K_8,K_9$ in such a way that  
 \begin{eqnarray}\nonumber
 \partial_t\big(\sum_0^i\partial_z\rho_t^{n,j}(\theta,u_{0,t}^{n,j}(\theta))\big)&=&\sum_0^i \partial_t\partial_z\rho_t^{n,j}(\theta,u_{0,t}^{n,j})+\partial_z^2\rho_t^{n,j}(\theta,u_{0,t}^{n,j})\partial_tu_{0,t}^{n,j}\\ \nonumber
\Big\|\partial_t\big(\sum_0^i\partial_z\rho_t^{n,j}(\theta,u_{0,t}^{n,j}(\theta))\big)\Big\|&\leq& T_{n,\infty}2C_{2}2^n\Gamma_{\alpha,\beta}(2^n)w_n +\frac{K_8}{\K2^n}\\\label{sirvex}
&\leq& \frac{K_9}{\K 2^n}.
 \end{eqnarray}
We will also need some estimates on the size of $\partial_t(\lambda_n^t)$. For this, it will be very useful to have a more developped expression for it. From  (\ref{valorlamda}) we have
 \begin{eqnarray*}
 \partial_t(\lambda_n^t)&=&\lambda_n^t\partial_t(\log \lambda_n^t)\\
 &=&\frac{\lambda_n^t}{\lambda}\Big(\int_{\T}\frac{\partial_t\rho_{1,t}^{n,i_n}}{1+\frac{\rho_{1,t}^{n,i_n}}{\lambda}}d\theta\Big)\\
 &=&\frac{\lambda_n^t}{\lambda}\Big(\int_{\T}\frac{\partial_t\rho_{1,t}^{n,0}+\partial_t\big(\sum_0^{i_n-1}\partial_z\rho_t^{n,j}(\theta,u_{0,t}^{n,j})\big)}{1+\frac{\rho_{1,t}^{i,i_n}}{\lambda}}d\theta\Big).
 \end{eqnarray*}
 By differentiating   (\ref{rhounonew}) we get 
\begin{equation*}
 \partial_t\rho_{1,t}^{n+1,0}= \partial_t(\lambda_n^t)e^{\{ \log (\rho_1^{n,i_n}+\lambda)\}|_{\overline{2^n}}}+\lambda_n^te^{\{ \log (\rho_{1,t}^{n,i_n}+\lambda)\}|_{\overline{2^n}}}\Big \{ \frac{\partial_t\rho_{1,t}^{n,i_n}}{\rho_{1,t}^{n,i_n}+\lambda}\Big \}\Big |_{\overline{2^n}}.
 \end{equation*}
The difference between the integral of $\partial_t\rho_{1,t}^{n+1,0}$ and   $\Delta_0$ may be written as the sum $I_1+I_2+I_3$ of three terms corresponding to the expressions in the three lines below: 
 \begin{eqnarray*}
  \int_{\T}\partial_t\rho_{1,t}^{n+1,0}(\theta)d\theta-\Delta_0 &=&\Big(\int_{T}\frac{\lambda_n^te^{\{ \log (\rho_1^{n,i_n}+\lambda)\}|_{\overline{2^n}}}}{\lambda}d\theta\Big)\Big(\int_{\T}  \frac{\partial_t\rho_{1,t}^{n,0}}{1+\frac{\rho_{1,t}^{n,i_n}}{\lambda}}d\theta  \Big)-\Delta_0\\
   &&+ \quad\Big(\int_{\T}\frac{\lambda_n^te^{\{ \log (\rho_1^{n,i_n}+\lambda)\}|_{\overline{2^n}}}}{\lambda}d\theta\Big) \Big(\int_{\T}  \frac{\partial_t\big(\sum_0^{i_n-1}\partial_z\rho_t^{n,j}(\theta,u_{0,t}^{n,j})\big)}{1+\frac{\rho_{1,t}^{n,i_n}}{\lambda}}d\theta  \Big)     \\
 &&+\quad\int_{\T}\lambda_n^te^{\{ \log (\rho_{1,t}^{n,i_n}+\lambda)\}|_{\overline{2^n}}}\Big \{ \frac{\partial_t\rho_{1,t}^{n,i_n}}{\rho_{1,t}^{n,i_n}+\lambda}\Big \}\Big |_{\overline{2^n}}d\theta .
 \end{eqnarray*}   
Notice that the  term $I_1$ itself equals the sum of the following three terms 
  \begin{eqnarray*}
 &&\Big(\int_{\T}\frac{\lambda_n^te^{\{ \log (\rho_1^{n,i_n}+\lambda)\}|_{\overline{2^n}}}}{\lambda}d\theta\Big)\Big(\int_{\T}  \frac{\partial_t\rho_{1,t}^{n,0}-1}{1+\frac{\rho_{1,t}^{n,i_n}}{\lambda}}d\theta  \Big)\\
 &&+\quad \Big(\int_{\T}\frac{\lambda_n^te^{\{ \log (\frac{\rho_1^{n,i_n}}{\lambda}+1)\}|_{\overline{2^n}}}}{\lambda}d\theta\Big)\Big(\int_{\T}  \frac{\frac{\rho_{1,t}^{n,i_n}}{\lambda}}{1+\frac{\rho_{1,t}^{n,i_n}}{\lambda}}d\theta  \Big)\\
 &&+\quad\Big(\int_{\T} e^{\int \log (\frac{\rho_1^{n,i_n}}{\lambda}+1)+\{ \log (\frac{\rho_1^{n,i_n}}{\lambda}+1)\}|_{\overline{2^n}}}-1 d\theta\Big)\Delta_0.
  \end{eqnarray*}
  Hence, there exists a constant $K_s$ such that
 \begin{equation}\label{KS}
 \|I_1\|\leq e^{K_{s}l_n}s_n+K_sl_n.
 \end{equation}
The terms  $\lambda_n^t,\rho_{1,t}^{n,i_n}$ are bounded from above by a universal constant. This together with (\ref{sirvex}) implies that there exists a constant   $C_6$ such that  
\begin{equation}
\|I_2\|\leq  \frac{C_6K_9}{\K 2^n}.
\end{equation}
The term $I_3$ equals the integral of a universally bounded term times the truncate of a function whose size is also bounded by a type $K$ constant. Notice that these bounds exist even before the introduction of any loss on the  width of the strip, see (\ref{silvio})). When introducing the loss on the width, we obtain a constant  $K_{10}$ for which the following bound holds: 
\begin{equation}
\|I_3\|\leq K_{10}\frac{e^{-2\pi 2^nd_n^1}}{d_n^1}=\frac{K_{10}l_{n+1}}{\K l_n}< \frac{K_{10}}{\K e^{2n+1}}.
\end{equation}
By considering these three bounds for the integrals  $I_1,I_2$ and $I_3$,  we get the following estimate at the third part of the stage: 
 \begin{equation}\label{vales2}
 \Big\|\int_{\T}\partial_t\rho_{1,t}^{n+1,0}d\theta-\Delta_0\Big\|\leq e^{K_sl_n}s_n+K_{s}l_n+\frac{C_6K_9}{\K 2^n}+\frac{K_{10}}{\K e^{2n+1}}.
    \end{equation}
    In the next paragraph we will give an explicit definition of the sequence  $\{s_n\}$ and also prove its claimed properties. 
  
  \paragraph{The sequence $\{s_n\}$.} We define this sequence by recurrence by letting
  \begin{eqnarray}
s_{n^*}&=&\frac{L^{-1}}{1000}\\
s_{n+1}&=&s_ne^{K_{s}l_n}+K_{s}l_n+\frac{C_6K_9}{\K 2^n}+\frac{K_{10}}{\K e^{2n+1}}.
\end{eqnarray}
Relation (\ref{vales2}) implies  (tautologically)  that 
\begin{equation}
\Big\|\int_{\T}\partial_t\rho_{1,t}^{n+1,0}d\theta-\Delta_0\Big\|\leq s_{n+1}.
\end{equation}
 By iterating the definition we get the estimate 
\begin{eqnarray*}
s_{n}&<& \Big(\prod_0^{\infty}e^{K_{s}l_j}\Big)s_{n^*}+\prod_0^{\infty}e^{K_{s}l_j}\Big (\sum_0^{\infty}K_{s}l_j+\frac{C_6K_9}{\K  2^n}+\frac{K_{10}}{\K e^{2j+1}}\Big )\\\
&\leq& \big(e^{K_{s}\sum_jl_j}\big )s_{n^*}+e^{K_{s}\sum_jl_j}\Big ( K_{s}\sum_0^{\infty}l_j+\frac{2C_6K_{9}}{\K}+\frac{K_{10}}{e\K}\Big ).
\end{eqnarray*} 
If $\K$ is large enough then this gives  (\ref{vales}) as desired $\quad_{\blacksquare}$
  \subsection{Fourth part, reduction on the parameter space}\label{coup}
  Now we arrive to a different stage of our process. Indeed, the estimates here do not arise neither as a consequence of the remainders of a linearized equation (quadratic residue) nor as the result of some Cauchy estimates associated to some loss in the strip. In this Section, we get our estimates as a consequence of the transversally properties of the family $\{F_t\}$. Let's see the necessity  of such a hypothesis.  From relation  (\ref{rhounonew}), at this moment we have a constant  $C_7$ such that
 \begin{equation}
\|\rho_{1,t}^{n+1,0}\| \leq C_{7}\|\rho_{1,t}^{n,i_n}\| \leq C_7K_4l_n.
\end{equation}
In other words, the  size of $\rho_{1,t}^{n+1,0}$ is still of order $l_n$. This fact is not a surprise, since we have  made operations which deal only with the size of the derivative   $\partial_{\theta}\rho_1$. We  will use the Affirmation  \ref{afirmac} to control $\rho_1$. Thus, we need to control the size of the complex number $\int_{\T}\rho_{1,t}(\theta)d\theta$. To do this we will strongly use the transversality property for the family  $\{F_t\}$. We will find a parameter $t_{n+1}$ such that $\int_{\T}\rho_{1,t_{n+1}}^{n+1,0}=0$. Then, we will exploit the fact that the size of  $\int_{\T}\rho_{1,t}^{n+1,0}$ is very small if we are placed very near this special parameter, due to the continuity of the functions. More precisely, as in (\ref{PN}) we will reduce the parameter space  $D(t_n,p_{n^-})$ by finding a new disc  $D(t_{n+1},p_{n+1})\subset D(t_n,p_{n^-})$ with $t_{n+1}$ being a simple zero of  $\int_{\T}\rho_{1,t}^{n+1,0}(\theta)d\theta$ and  $p_{n+1}$ small enough.
\subsubsection{An application of Rouch\'e Lemma}
We will use the classical Rouch\'e Lemma to find $t_{n+1}$. To do this, we need to show that the hypothesis of the Lemma are satisfied on the boundary of some disc around  $t_n$.  We develop
\begin{eqnarray*}
\int_{\T}\frac{\rho_{1,t}^{n+1,0}(\theta)}{\lambda}d\theta&=&\int_{\T}\Big(\frac{\lambda_n^t}{\lambda} e^{ \{\log (\rho_{1,t}^{n,i_n}(\theta)+\lambda) \}|_{\overline{2^n}}}-1\Big)d\theta\\
&=&\int_{\T}\Big(e^{\int_{\T} \log \big (\frac{\rho_{1,t}^{n,i_n}(\theta)}{\lambda}+1 \big )d\theta + \{ \log(\rho_{1,t}^{n,i_n}(\theta)+\lambda )\}|_{\overline{2^n}}}-1\Big)d\theta.
\end{eqnarray*}
 We then define  $g_+(t)$ as being the value of this last expression. We compare  this function with the function $g(t)=\int_{\T}\frac{\rho_{1,t}^{n,0}(\theta)}{\lambda}d\theta$ on the boundary of some disc  $D(t_n,r)$:
 \begin{eqnarray*}
 |g_+(t)-g(t)|&\leq&\int_{\T}\bigg\{\Big|\int_{\T}\log \Big(\frac{\rho_{1,t}^{n,i_n}(\theta)}{\lambda}+1\Big)d\theta-g(t)
 + \big \{\log (\rho_{1,t}^{n,i_n}(\theta)+\lambda)\big \}|_{\overline{2^n}}\Big| +Z^2e^{Z}\bigg\}d\theta,
 \end{eqnarray*}
where $Z=\big|\int_{\T}\log \big(\frac{\rho_{1,t}^{n,i_n}(\theta)}{\lambda}+1\big)d\theta+\big \{\log (\rho_{1,t}^{n,i_n}(\theta)+\lambda)\big \}|_{\overline{2^n}}\big|$. Then there exists a constant  $K_{11}$  such that  $Z<K_{11}l_n<1$,  provided that $\K$ is large enough. Thus, using the Taylor series expansion of  $\log\Big (\frac{\rho_{1,t}^{n,i_n}(\theta)}{\lambda}+1\Big)$, we have
  \begin{eqnarray*}
  |g_+(t)-g(t)|&\leq&\Big |\int_{\T} \Big(\frac{\rho_{1,t}^{n,i_n}(\theta)}{\lambda}-\frac{\omega(\theta,t)}{2}\Big(\frac{\rho_{1,t}^{n,i_n}(\theta)}{\lambda}\Big )^2-\frac{\rho_{1,t}^{n,0}(\theta)}{\lambda}\Big)d\theta\Big|\\
 &&+\int_{\T}\Big|\Big\{\log\Big( \rho_{1,t}^{n,i_n}(\theta)+\lambda\Big)\Big \}\Big|_{\overline{2^n}}+Z^2e^Z\Big |d\theta\\
&<& \int_{\T}\Big|\sum_0^{i_n-1}\partial_z\rho_t^{n,j}(\theta,u_{0,t}^{n,j}(\theta))\Big|+\Big| \frac{\omega(\theta,t)}{2}\Big(\frac{\rho_{1,t}^{n,i_n}(\theta)}{\lambda}\Big )^2\Big|d\theta\\
&&+\int_{\T}\Big|\Big\{\log\Big( \rho_{1,t}^{n,i_n}(\theta)+\lambda\Big)\Big \}\Big|_{\overline{2^n}}+Z^2e^Z\Big |d\theta,
 \end{eqnarray*}
where the absolute value of  $\omega(\theta,t)$ lies between $0$ and $1$.
One can easily check the following estimates:
\begin{eqnarray}\label{tula123}
\int_{\T}\Big|\sum_0^{i_n-1}\partial_z\rho_t^{n,j}(\theta,u_{0,t}^{n,j}(\theta))\Big|d\theta&\leq&\frac{2C_2K_Ml_n}{\K}\\\label{tula2}
\int_{\T} \Big|\frac{\omega(\theta,t)}{2}\Big(\frac{\rho_{1,t}^{n,i_n}(\theta)}{\lambda}\Big )^2\Big|d\theta&\leq&\frac{K_4^2l_n^2}{2}\\\label{tula3}
\int_{\T}\Big|\Big\{\log\Big( \rho_{1,t}^{n,i_n}(\theta)+\lambda\Big)\Big \}\Big|_{\overline{2^n}}\Big| d\theta &\leq& \frac{C_3K_4l_ne^{-2\pi2^nd_n^1}}{d_n^1}\leq \frac{C_3K_4l_{n+1}}{\K}.
\end{eqnarray} 
Therefore, if $\K$ is large enough,  for every  $r\leq p_{n^-}$ 
and every $t \in D(t_n,r)$ we have 
\begin{equation}\label{tula4}
|g_+(t)-g(t)|<3l_nL^{-1}.
\end{equation}
The transversality condition on the derivative  $\partial_t\int_{T}\rho_{1,t}(\theta)d\theta$ implies that $g(t)$ grows almost as a linear application, namely, 
 \begin{equation*}
 s_n|t-t_n|\geq\Big | \int_{t_n}^t\big(\frac{\Delta_0}{\lambda}-\partial_tg(t)\big)\Big |=\Big|\frac{\Delta_0(t-t_n)}{\lambda}-g(t)+g(t_n)\Big|.
 \end{equation*}
This implies that the inequality 
 \begin{equation*}
 |g(t)|\geq\big (L^{-1}-s_n\big )|t-t_n|\geq \frac{99}{100}L^{-1}|t-t_n|=\frac{99}{100}L^{-1}r
 \end{equation*} 
holds on the boundary of a disc $D(t_n,r)$. Therefore, if the radius  $r$ is at least   $\frac{100}{99}3l_n$, then the hypothesis in the Rouch\'e Lemma is satisfied. So, we get a simple zero,   $t_{n+1}$, for  $g_+(t)$ inside $D(t_n,r)$.  Now we need to show that if we choose $r$ as being this minimal value, then the new zero $t_{n+1}$ lies effectively inside the parameter space $D(t_n,p_{n^-})$. To do this, it is enough to see that 
\begin{equation}\label{tula1}
p_{n^-}-3l_n\frac{100}{99}>92l_n>0.
\end{equation}
Moreover, we need to show that we can start with the stage $n+1$ having a parameter space disc of radius  $p_{n+1}$. For this we notice that
\begin{eqnarray}\nonumber
p_{n+1}&=&100l_{n+1}\\\nonumber
&=&l_n\Big(\frac{100\Gamma_{\alpha,\beta}(2^n)\Gamma_{\alpha}(2^n)}{\Gamma_{\alpha,\beta}(2^{n+1})\Gamma_{\alpha}(2^{n+1})2e^{2n+1}}\Big)\\\label{yapoh}
&<&25l_n.
\end{eqnarray}
From  (\ref{tula1}) and  (\ref{yapoh}) we conclude that the disc  $D(t_{n+1},p_{n+1})$ is contained in $D(t_n,p_{n^-})$. Thus, the new parameter space is  the disc  $D(t_{n+1},p_{n+1})$. 
\newline
 In this way we get  the whole necessary estimates in order to relaunch our process at the stage $n+1$, and this finishes the stage $n$.
\subsection{Convergence of the method}\label{convergencia} 
As we have seen in the previous Sections, we can iterate the algorithm infinitely times for the \emph{fhd} corresponding to the parameter $\bar{t}=\cap_{n\geq n^*}D(t_n,p_n)$. (Of course, $\bar{t}$ is the parameter claimed  in the Theorem \ref{teo}.) In what follows, we will omit this parameter in all of the notation, and we will deal with 
$F = F_{\bar{t}}$. We need  to show that the successive compositions of the conjugacies constructed in the  course of our process gives raise to an analytical change of coordinates. The decreasing rates  for the functions and all of the estimates arising in our process  are sufficient for this. We let 
\begin{equation}
\delta^{\infty}=\delta^*-\sum_{n^*}^{\infty}\max(d_n^0,d_n^1)>\frac{\delta}{2},
\end{equation}
and we  denote  
\begin{equation*}
H_n(\theta,z)=h_{n,0}\circ h_{n,1}\circ\dots \circ h_{n,i_n-1}\circ h_n(\theta,z)=\Big(\theta,u_0^n(\theta)+u_1^n(\theta)z\Big),
\end{equation*}
where we write   $u_0^n(\theta)=\sum_{i=0}^{i_n-1}u_0^{n,i}(\theta)$. The successive compositions of the conjugacies up to the stage $n$ is 
\begin{eqnarray*}
H_{n^*}\circ H_{n^*+1} \circ \dots \circ H_n(\theta,z)&=&\big (\theta,u_0^{n^*}+u_1^{n^*}\big ( u_0^{n^*+1}+u_1^{n^*+1}(\ldots(u_0^n+u_1^nz)\ldots)\big)\big)\\
&=&	\big (\theta, u_0^{n^*}+u_1^{n^*}u_0^{n^*+1}+\ldots+\left (\prod_{j=n^*}^{n-1}u_1^j\right )u_0^n+\left( \prod_{j=n^*}^nu_1^j\right )z\big).
\end{eqnarray*}

We need to show that the constant term and the coefficient of $z$ above,  converge uniformly on $\theta$ when $n$ goes to infinity. For this, it is enough to see that the product  $\prod_{n^*}^{\infty}u_1^j(\theta) $ converges uniformly and the  series  $\sum_{n^*}^{\infty}u_0^j(\theta)$ converges absolutely and uniformly. For the product we have
\begin{equation*}
\log\left( \prod_{j=n^*}^iu_1^j(\theta)\right )= \sum_{j=n^*}^{i}v^j(\theta).
\end{equation*}
The above series converges absolutely and uniformly by the estimate 
 (\ref{vvv}) and point 2. of the Lemma \ref{ele1}. The series of  $u_0$ is
\begin{equation*}
\sum_{j=n^*}^n u_0^j(\theta)=\sum_{j=n^*}^n\sum_{i=0}^{i_j-1}u_0^{j,i}(\theta),
\end{equation*}
which allows to obtain the estimate
\begin{equation}\label{tallau}
\Big \|\sum_{j=n^*}^n u_0^j\Big \|_{B_{\delta^{\infty}}}\leq\sum_{j=n^*}^n\sum_{i=0}^{i_j-1}C_22^j\Gamma_{\alpha,\beta}(2^j)\frac{w_j}{2^i}<\frac{1}{8}\sum_{j=n^*}^{\infty}\frac{1}{4^je^{j^2}}.
\end{equation}
  Thus, the following limit exists 
 \begin{equation*}
 H(\theta,z)=\lim_{n\to \infty}H_{n^*}\circ H_{n^*+1} \circ \dots \circ H_n(\theta,z).
 \end{equation*}
We write this limit in the form
 \begin{equation*}
 H(\theta,z)=(\theta,u_0(\theta)+u_1(\theta)z),
 \end{equation*}
 where the functions  $u_0,u_1$ are analytical on the strip  $B_{\delta^{\infty}}$. The function  $u_1$ has zero topological degree since the topological degree is a homomorphism and each function   $u_1^j$ has zero degree. The conjugacy of  $F$ by $H$  takes the  form 
 \begin{equation*}
 H^{-1}\circ F\circ H(\theta,z)=(\theta+\alpha,\lambda z + \rho^{\infty}(\theta,z)).
 \end{equation*}
  Indeed, the sequences  $l_n,w_n$ which control the size of   $\rho_0^{n,0}, \rho_1^{n,0}$ go to  $0$ when $n$ goes to infinity. The function  $\rho^{\infty}$ vanishes up to the order  2 at $z=0$ and it is defined on the disc  $D(0,R^{\infty})$,  with $\frac{3}{8}\leq R^{\infty}$. The invariant curve is  $u_0(\theta)$, and due to (\ref{tallau}) it verifies 
  \begin{equation}\label{talladeufinal}
  \|u_0\|_{B_{\delta^{\infty}}}\leq \frac{1}{8}\sum_{j=n^*}^{\infty}\frac{1}{4^je^{j^2}}.
  \end{equation}
    \subsection{A modification at the first stage and the size of  $\eps$}\label{remarcaepsilon}
    Let  $\bar{n}\geq n^*$ be a natural number and $\eps$ be a real positive number such that 
   \begin{equation}
    \frac{\delta}{24}l_{\bar{n}+1}<\eps\leq \frac{\delta}{24}l_{\bar{n}}.
    \end{equation}
    Let  $\{F_t\}_{t\in \Sigma}$ be an analytical family verifying the hypothesis of our Main Theorem and 
    \begin{equation}\label{anal199}
    \|\rho_1\|\leq \eps\quad ,\quad \|\rho_0\|\leq \frac{24\eps}{\delta \K 4^{\bar{n}}\gab(2^{\bar{n}})}.
    \end{equation}  
    We start the algorithm at the stage    $\bar{n}$. We notice   that  in order to pass to the stage $\bar{n}+1$ we need to show the following facts at the end of the stage $\bar{n}$:
     \paragraph{In the first part. } The inequalities  (\ref{iW}), (\ref{iL}) must be 
     \begin{eqnarray}
     \|\eta_{\bar{n},i}^t\|&\leq&\frac{24\eps}{\delta\K 4^{\bar{n}}\gab(2^{\bar{n}})2^i}\\
     \|\partial_{\theta}\rho_{1,t}^{\bar{n},i}\|&\leq&2\big(\frac{24\eps}{\delta}\big).
     \end{eqnarray}
     \paragraph{In the second part. } Inequality (\ref{Kocho}) must be 
     \begin{equation}
     \|u_{1,t}^{\bar{n}}-1\|\leq C_3C_2C_42^{\bar{n}}\ga(2^{\bar{n}})K_4\frac{24\eps}{\delta}.
     \end{equation}
      Now the first two inequalities follow by repeating  adequately  the computations of the first part, while the third one follows immediately.
    \\
     
     In this way, we reach   the third part and we obtain successfully the corresponding bounds for 
     $l_{\bar{n}+1}$ and $w_{\bar{n}+1}$ when introducing the loss in the strip.
     \paragraph{In the estimates with respect to the parameter. } In order to get (\ref{nosequewa}) it is enough at each step to alow a loss of size  $\frac{24\eps}{\delta}\Big(\frac{2}{3}\Big)^i$ for the radius of the parameter space disc. Moreover, to obtain (\ref{silvio}) it is enough to alow a loss of size  $\frac{24\eps}{\delta}$. At the fourth part we can see that  (\ref{tula123}), (\ref{tula2}), (\ref{tula3}), are bounded from above by  
     \begin{equation}
     \frac{2C_2K_M24\eps}{\K \delta}\quad, \quad \frac{K_4^2\big(\frac{24\eps}{\delta}\big)^2}{2}\quad, \quad \frac{C_3K_4\big(\frac{24\eps}{\delta}\big)e^{-2\pi 2^{\bar{n}}d_{\bar{n}}^1}}{d_{\bar{n}}^1}\leq \frac{C_3K_4l_{\bar{n}+1}}{\K} 
     \end{equation} 
     respectively. If $\K$ is large enough  then (\ref{tula4}) is bounded from above by $3L^{-1}\Big(\frac{24\eps}{\delta}\Big)$. This allows us to finish the stage  with an appropriate loss of width for the strip. 
     \\
     
     The above considerations show that in order to start with the stage $\bar{n}$ we just  need a parameter space disc of radius 
     \begin{equation}
     p_{\bar{n}}=100\frac{24\eps}{\delta}.
     \end{equation}    
For such a disc, we also verify that the parameter $\bar{t}$ proving  the existence of the invariant curve is located  at a distance  $\eps$ times a constant from the origin in the complex plane, as desired. If $\eps$ goes to  $0$, then the level  $\bar{n}$ of the first stage grows and   (\ref{talladeufinal}) allows to conclude that the size of the invariant curve tends to  $0$.
 \\    
    
Actually, we have shown that Theorem  \ref{teo} holds under the stronger hypothesis  $\|\rho_0\|\leq w_{n^*}$ (or more precisely under hypothesis  (\ref{anal199})). In the next Section we will see that we only need to impose the hypothesis  $\|\rho_0\|\leq \eps$. We will see that  even under this weaker hypothesis  we can reduce the problem to  the stronger hypothesis used up so far.    Indeed, we will describe a \emph{previous} preparative process  showing this reduction.
\section{A previous preparative process and the end of the proof}\label{previous}
    \subsection{A initial dicrease in the size of  $\rho_0$}\label{epsiloncuadrado}
 We know that solving the $u_0$ equation gives rise to a diminution on the size of $\rho_0$. Thus, we will solve many times this equation  to get the adequate size for  $\rho_0$  in order to start the algorithm. Although each time we solve the equation  we slightly increase  the size of  $\rho_1$, the total amount of the increase is well behaved (and may be effectively controlled). We will omit many details when manipulating these estimates, the corresponding justifications being contained in the previous  Sections. 
 
 We fix a constant  $\K$ by using the data  $\big(L,M,2T+(2000L)^{-1},\delta,(\alpha,\beta)\big)$ as described  at the begining of the Section  \ref{aleph}. This choice gives us 
our sequences $\{l_n\},\{w_n\},\{d_n^0\},\{d_n^1\}$ and a natural number  $n^*$ representing the lowest starting stage for the algorithm.  Suppose that there exists  $\bar{n}\geq n^*$ and $\eps>0$ such that 
 \begin{equation}
 \frac{\delta l_{\bar{n}+1}}{24}<\eps (w_{n^*})^{-1}\leq \frac{\delta l_{\bar{n}}}{24}.
 \end{equation}
Assume also that the hypothesis of Theorem  \ref{teo} hold and that  for every parameter  $t$ in the disc  $D(0,K_R\eps)$  we have 
 \begin{equation}
 \|\rho_{0,t}\|\leq \eps\quad, \quad \|\rho_{1,t}\|\leq \eps,
 \end{equation}
where $K_R$ is a constant which will be defined in a sequel. The value of the radius for  the parameter space disc corresponds exactly to the sum of what we need to start this previous preparative process and the necessary radius to start with the algorithm at the stage  $\bar{n}$.  This previous preparative process also has stages (indexed by the natural numbers between $n^*$ and $\bar{n}-1$). Each one of these stages contains two parts. The first part of the stage $n$ is at some moment divided into steps. Each step consists in solving  the equation of $u_0$ once, with truncation up to the order  $2^n$. The second part  consists in introducing  a loss for the width of the strip of size  $d_n^0$. For notation reasons we will use the sequence of real numbers
 \begin{equation}
 W_n=C_2M2^n\gab(2^n).
 \end{equation}
 We present  the first stage $n^*$ separately. The reason for doing this is twofold. On the one hand, at this stage there the estimates are slightly different. On the other hand, this stage allows us to  understand the manipulations and estimates that arise in the remaining stages. As usual (in this work), we denote the starting functions as   $\rho_{0,t}^{n^*},\rho_{1,t}^{n^*},\rho_t^{n^*}$. We introduce an index to indicate the steps. We define the sequence of functions    $\eta_{n^*,i}^t$ by
 \begin{eqnarray}
 \eta_{n^*,0}^t&=&\rho_{0,t}^{n^*,0}\\
 \eta_{n^*,i+1}^t&=&\rho_{1,t}^{n^*,i}u_{0,t}^{n^*,i}+\rho_t^{n^*,i}(\cdot,u_{0, t}^{n^*,i}).
 \end{eqnarray}
 The hypothesis of Theorem \ref{teo} say that $\|\eta_{n^*,0}^t\|\leq \eps$. We will show inductively that 
 \begin{eqnarray}
 \|\eta_{n^*,i}^t\|&\leq&\frac{\|\eta_{n^*,0}\|}{2^i}\leq \frac{\eps}{2^i}\\
 \|\rho_{1,t}^{n^*,0}\|&\leq& \|\rho_{1,t}^{n^*,0}\|+2\|\eta_{n^*,0}^t\|W_{n^*}\leq \eps(1+2W_{n^*})< \eps(w_{n*})^{-1}.
 \end{eqnarray}
To do this, assume that these inequalities hold for  $0\leq j<i+1$. Then we can estimate the new  $\eta$ by
 \begin{equation*}
 \|\eta_{n^*,i+1}^t\|\leq \eps(1+2W_{n^*})\|\eta_{n^*,i}^t\|\frac{W_{n^*}}{M}+M\Big(\|\eta_{n^*,i}^t\|\frac{W_{n^*}}{M}\Big)^2.
 \end{equation*} 
 Since $\eps(1+2W_{n^*})<\eps(w_{n^*})^{-1}, \|\eta_{n^*,i}^t\|\leq \eps$ and  $\bar{n}\geq n^*$ we have
 \begin{equation*}
 \eps(1+2W_{n^*})\frac{W_{n^*}}{M}<l_{\bar{n}}\frac{\delta W_{n^*}}{24M}<\frac{1}{4}\quad,\quad\|\eta_{n^*,i}^t\|\frac{W_{n^*}^2}{M}\leq \eps\frac{W_{n^*}^2}{M}<\frac{1}{4}.
 \end{equation*}
 This implies that  $\|\eta_{n^*,i+1}^t\|\leq \|\eta_{n^*,i}\|2^{-1}$. The size of  $\|\rho_{1,t}^{n^*,i+1}\|$ is bounded from above by
 \begin{equation*}
 \|\rho_{1,t}^{n^*,i+1}\|\leq \|\rho_{1,t}^{n^*,0}\|+W_{n^*}\sum_0^i\|\eta_{n^*,j}^t\|< \|\rho_{1,t}^{n^*,0}\|+2W_{n^*}\|\eta_{n^*,0}^t\| \leq \eps(1+2W_{n^*}).
 \end{equation*}
  We continue the steps until to  the moment when the size of $\eta$ is  smaller than  $\frac{\eps w_{n^*+1}}{2\K w_{n^*}}$, or more precisely when   $\frac{\eps}{2^i}\leq \frac{\eps w_{n^*+1}}{2\K w_{n^*}}$ (this subtle difference will be very important in Section  \ref{teopara}). Let $i_{n^*}$ be such a moment. In the second part of the stage we introduce a loss in the width of the strip of size  $d_{n^*}^0$. The final size of  $\rho_0$ is thus bounded by
 \begin{eqnarray*}
 \|\rho_{0,t}^{n^*+1}\|&\leq& \|\eta^{t}_{n^*,i_{n^*}}\|+2\|\eta_{n^*,0}^t\|\frac{e^{-2\pi 2^{n^*}d_{n^*}^0}}{d_{n^*}^0}\\
 &\leq& \frac{\eps w_{n^*+1}}{2\K w_{n^*}}+\frac{2\eps w_{n^*+1}}{4\K w_{n^*}} \\
 &\leq&\frac{\eps w_{n^*+1}}{\K w_{n^*}}.
 \end{eqnarray*}
 \subsubsection{An arbitrary  stage  $n^*<n\leq \bar{n}-1$} 
 At the begin of the stage we have 
 \begin{eqnarray}
 \|\rho_{0,t}^{n,0}\|&\leq&\frac{\eps w_n}{\K w_{n^*}}\\
 \|\rho_{1,t}^{n,0}\|&\leq& \eps\Big (1+2W_{n^*}+\frac{2w_{n^*+1}W_{n^*+1}}{\K w_{n^*}}+\dots + \frac{2w_{n-1}W_{n-1}}{\K w_{n^*}}\Big)\\&<&\eps(w_{n^*})^{-1}.
 \end{eqnarray} 
 For the usual sequence $\{\eta_{n,i}^t\}$ we will show inductively that  
 \begin{eqnarray}
 \|\eta_{n,i}^t\|&\leq& \frac{\|\eta_{n,0}^t\|}{2^i}\\
 \|\rho_{1,t}^{n,i}\|&\leq&\eps\Big(1+2W_{n^*}+\dots+\frac{2w_nW_n}{\K w_{n^*}} \Big)<\eps(w_{n^*})^{-1}. 
 \end{eqnarray}
 So, suppose these inequalities hold for $0\leq j< i+1$. The size  of $\eta$ is bounded from above by
 \begin{equation*}
 \eps(w_{n^*})^{-1}\frac{W_n}{M}\|\eta_{n,i}^t\|+M\Big(\|\eta_{n,i}^t\|\frac{W_n}{M}\Big)^2.
 \end{equation*}
 Moreover, we have 
 \begin{eqnarray*}
 \eps(w_{n^*})^{-1}\frac{W_n}{M}&<& l_{n}\frac{\delta W_n}{24 M}<\frac{1}{4}\\
 M\|\eta_{n,i}^t\|W_n^2&\leq& \frac{M\eps w_{n}W_n^2}{\K w_{n^*}}<\frac{\delta M}{24 \K} l_nw_nW_n^2<\frac{1}{4}. 
 \end{eqnarray*}
 That implies  $\|\eta_{n,i+1}^t\|\leq \|\eta_{n,i}^t\|2^{-1}$. The size of  $\rho_{1,t}^{n,i+1}$ is bounded  by
 \begin{eqnarray*}
 \|\rho_{1,t}^{n,i+1}\|&\leq& \|\rho_{1,t}^{n,0}\|+W_n\sum_0^i\|\eta_{n,j}^t\|\\
 &<& \|\rho_{1,t}^{n,0}\|+\frac{2W_n\eps w_n}{\K w_{n^*}}\\
 &\leq& \eps\Big(1+2W_{n^*}+\dots+ \frac{2W_nw_n}{\K w_{n^*}}\Big)<\eps (w_{n^*}){-1}.
 \end{eqnarray*}
 We continue the steps up to the first moment where the size of  $\eta$ is  smaller than  $\frac{\eps w_{n+1}}{2\K w_{n^*}}$, or more precisely when  $\frac{\eps w_n}{\K w_{n^*}2^i}\leq \frac{\eps w_{n+1}}{2\K w_{n^*}}$. Let  $i_n$ be such moment. At the second part of the stage we introduce a loss in the width of the strip of size  $d_n^0$. The final size of  $\rho_0$ is bounded by
 \begin{eqnarray*}
 \|\rho_{0,t}^{n+1}\|&\leq& \|\eta^{t}_{n,i_{n}}\|+2\|\eta_{n,0}^t\|\frac{e^{-2\pi 2^{n}d_{n}^0}}{d_{n}^0}\\
 &\leq& \frac{\eps w_{n+1}}{2\K w_{n^*}}+2\frac{\eps w_n}{\K w_{n^*}}\frac{\eps w_{n+1}}{4\K w_{n}} \\
 &<&\frac{\eps w_{n+1}}{\K w_{n^*}}.
 \end{eqnarray*}
 Hence, at the end of the  $(\bar{n}-1)$-th stage we obtain  
 \begin{eqnarray}
 \|\rho_{0,t}^{\bar{n}}\|&\leq& \frac{\eps w_{\bar{n}}}{\K w_{n^*}}< \frac{\delta \eps (w_{n^*})^{-1}}{\K 4^{\bar{n}}\gab(2^{\bar{n}})}\\
 \|\rho_{1,t}^{\bar{n}}\|&\leq& \eps (w_{n^*})^{-1}.
 \end{eqnarray}
 These functions are defined on the strip  $\delta-\sum_{n^*}^{\bar{n}-1}d_n^0$. The function $\rho_t^{\bar{n}}(\theta,\cdot)$ is defined, for the $z$ variable, on a disc of radius 
 \begin{equation}
 1-\sum_{n=n^*}^{\bar{n}-1}\sum_{i=0}^{i_n}\|u_{0,t}^{n,i}\|.
 \end{equation} 
 If we start the algorithm at the stage  $\bar{n}$ with this parameter space disc, we will reach  the end of the process with a radius even larger than  $3/8$. The size of  $\|\rho_t\|$ does not change with the operation of resolution of the equation of  $u_0$.
 \subsubsection{Estimates with respect to the parameter}
 We will need to estimate the values of the following two sums: 
     \begin{equation*}
     \sum_{n=n^*}^{\bar{n}-1}\sum_{i=0}^{i_n}\|u_{0,t}^{n,i}\|\quad ,\quad \sum_{n=n^*}^{\bar{n}-1}\sum_{i=0}^{i_n}\|\partial_t u_{0,t}^{n,i}\|.
     \end{equation*}
     The first one is bounded by
     \begin{eqnarray*}
      \sum_{n=n^*}^{\bar{n}-1}\sum_{i=0}^{i_n}\|u_{0,t}^{n,i}\|&\leq&  2\frac{W_{n^*}\eps}{M}+2\sum_{ n>n^*}\frac{W_n\eps w_n}{M\K w_{n^*}}\\
      &\leq&\frac{2W_{n^*}\eps}{M}\Big( 1+\frac{1}{2\K}\sum_{n>n^*}\frac{l_n}{l_{n^*}}\Big)<\frac{4W_{n^*}}{M}\eps.
           \end{eqnarray*}
            By introducing a loss of size  $96000LW_{n^*}\eps$ on the radius of the parameter space disc,  we obtain a Cauchy estimate
           \begin{equation}
             \sum_{n=n^*}^{\bar{n}-1}\sum_{i=0}^{i_n}\|\partial_t u_{0,t}^{n,i}\|\leq \frac{1}{4000LM}.
           \end{equation}
           Now, we can estimate the final oscillation of  $\partial_t\partial_z \rho_t$ by
           \begin{equation}
           osc\big(\partial_t\partial_z\rho_t^{\bar{n}}\big)\leq 2T +2M   \sum_{n=n^*}^{\bar{n}-1}\sum_{i=0}^{i_n}\|u_{0,t}^{n,i}\|\leq 2T+\frac{1}{2000L}.
           \end{equation} 
           We will need the following estimates
                      \begin{eqnarray*}
           \Big\|\partial_t\Big(\sum_{n=n^*}^{\bar{n}-1}\sum_{i=0}^{i_n}\partial_z\rho_t^{n,i}(\cdot,u_{0,t}^{n,i})\Big)\Big\|&\leq& \Big(2T+\frac{1}{2L}\Big)\sum\sum\|u_0\|+M\sum\sum\|\partial_tu_0\|\\
           &\leq& \Big(2T+\frac{1}{2L}\Big)\frac{4W_{n^*}}{M}\eps+\frac{1}{4000L}\\
           &\leq&\frac{1}{2000L}
           \end{eqnarray*}
           where the last inequality follows since the size of $\K$ is very large.
           We define the complex number  $\Delta_0$ by 
                      \begin{equation}
           \Delta_0=\partial_t\Big(\int_{\T}\rho^{n^*}_{1,t}d\theta \Big)\Big|_{t=t_0}.
           \end{equation}
  \subsection{A simple zero for  $\int_{\T}\rho_{1,t}$}\label{zerosimple}
  There is  only one remaining ingredient we miss in order to apply the algorithm to our family, namely, the existence of a simple zero for  $\int_{\T}\rho_{1,t}^{\bar{n}}$ at the center of the parameter space disc.
We will find this zero by using Rouch\'e Lemma. Thus, we need to compare  $\int_{\T}\rho_{1,t}^{\bar{n}}$ with the linear part of  $\int_{\T}\rho_{1,t}^{n^*}$ at $t_0$ on the boundary of a disc of radius  $R$ and centered at  $t_0$:  
  \begin{eqnarray*}
  \Big|\int_{\T}\rho_{1,t}^{\bar{n}}-\partial_t\Big(\int_{\T}\rho_{1,t}^{n^*}\Big)\Big|_{t=t_0}(t-t_0)\Big|&\leq&\Big|\int_{\T}\rho_{1,t}^{n^*}-\partial_t\Big(\int_{\T}\rho_{1,t}\Big)\Big|_{t=t_0}(t-t_0)\Big|\\
  &&+\Big|\int_{\T}\rho_{1,t}^{\bar{n}}-\int_{\T}\rho_{1,t}^{n^*}\Big|\\
   &\leq&\eps+\frac{T}{2}R^2+4W_{n^*}\eps.
  \end{eqnarray*}
We have the linear estimate  
  \begin{equation*}
  \Big|\partial_t\Big(\int_{\T}\rho_{1,t}^{n^*}\Big)\Big|_{t=t_0}(t-t_0)\Big|>L^{-1}R
  \end{equation*}
  over this boundary. We choose  the radius   
  \begin{equation}
  R_{\eps}=\frac{L^{-1}-\sqrt{L^{-2}-2T\eps(1+4W_{n^*})}}{T},
  \end{equation} 
  which corresponds to the smallest root of the equation between the two precedent bounds. In this way we assure the existence of a simple zero  $t^*$  of $\int_{\T}\rho_{1,t}^{\bar{n}}$ inside $D(t_0,R_{\eps})$.   We easily verify that  $R_{\eps}$ is a positive real number. Furthermore,  we also verify that
  \begin{equation}
  R_{\eps}<\frac{L^{-1}}{T}\frac{2T\eps(1+4W_{n^*})}{L^{-2}}=2L\eps(1+4W_{n^*}),  
  \end{equation} 
  because  $x>1-\sqrt{1-x}$ for every small $x$. 
  \\
  
  In order to start the iterative part of our process (the algorithm described in Sections \ref{realisacion}$\dots$\ref{coup}) we need a parameter space disc of radius 
     $100\frac{24\eps(w_{n^*})^{-1}}{\delta}$ and centered at  $t^*$. All what we have seen allows us to define the constant  $K_R$ as being precisely what is needed. That is,   to find   $t^*$, and   to lose a  little bit of radius to obtain the estimates with respect to the parameter in the previous preparative process. Thus we have to choose 
  \begin{equation}
K_R=2L(1+4W_{n^*})+96000LW_{n^*}+\frac{2400(w_{n^*})^{-1}}{\delta}.
  \end{equation}
 To finish this previous preparative process and to start with the algorithm, it is enough to notice that the following estimate
  \begin{eqnarray*}
  \Big\|\Delta_0-\partial_t\Big(\int_{\T}\rho_{1,t}^{\bar{n}}\Big)\Big\|&\leq& \Big\|\Delta_0-\partial_t\Big(\int_{\T}\rho_{1,t}^{{n^*}}\Big)\Big\|\\
  &&+ \Big\|\partial_t\Big(\sum_{n=n^*}^{\bar{n}-1}\sum_{i=0}^{i_n}\partial_z\rho_t^{n,i}(\cdot,u_{0,t}^{n,i})\Big)\Big\|\\
  &\leq&T\eps\Big(2L(1+4W_{n^*})+\frac{2400(w_{n^*})^{-1}}{\delta}\Big) +\frac{L^{-1}}{2000}\\
  &\leq&\frac{L^{-1}}{1000}
  \end{eqnarray*}
  holds inside the parameter space disc centered at $t_0$ with radious  $\eps\Big(2L(1+4W_{n^*})+\frac{2400(w_{n^*})^{-1}}{\delta}\Big)$. Notice that this disc contains the disc centered at  $t^*$ with radius  $\eps\frac{2400(w_{n^*})^{-1}}{\delta}$.
  
 The invariant curve is  
  \begin{equation}\label{formuladeucero}
  u=\sum_{n=n^*}^{\bar{n}-1}\sum_{i=0}^{i_n}u_{0,\bar{t}}^{n,i}+\lim_{\tilde{n}\to \infty}\sum_{n=\bar{n}}^{\tilde{n}}u_{0,\bar{t}}^n\Big(\prod_{j=\bar{n}}^n u_{1,\bar{t}}^j\Big).
  \end{equation}
  It has zero degree because all the coordinate changes performed in this previous preparative process have zero degree $\quad_{\blacksquare}$ 
\section{A parametrized version of the persistence of the invariant curve Theorem for the fibered holomorphic dynamics}\label{teopara}
Let $\Lambda\in \C$ be an open set. We consider an application from  $\Lambda$ to the set of analytical $1$ complex parameter families of fibered holomorphic dynamics,
\begin{equation*}
s\in \Lambda \longmapsto \{F^s_t\}_{t\in \Sigma}.
\end{equation*} 
We say that this application is analytic if 
\begin{equation*}
(s,t,\theta,z)\longmapsto F^s_t(\theta,z)
\end{equation*}
is an analytical function. Let  $\bar{s}$ be in   $\Lambda$. Suppose that the family  $\{F_t^{\bar{s}}\}_{t\in \Sigma}$ verify the hypothesis of  Theorem \ref{teo}.  Hence, for every  $t$ in $D(t_0,K_R\eps)\subset \Sigma$ one has
\begin{itemize}
\item[$\bullet$] $L>\Big|\partial_t\Big(\int_{\T}\rho_{1,t}^{\bar{s}}\Big)\big|_{t=t_0}\Big|>L^{-1}$
\item[$\bullet$] $\|\rho_{0,t}^{\bar{s}}\|< \eps\quad,\quad \|\rho_{1,t}^{\bar{s}}\|< \eps$
\item[$\bullet$] $\|\partial_z^2\rho_t^{\bar{s}}\|<M\quad ,\quad \|\partial_t\partial_z\rho_t^{\bar{s}}\|+\|\partial_t^2\rho_{1,t}^{\bar{s}}\|<T$
\end{itemize}
for some  $\eps\in(0,\eps^*]$, where  $\eps^*$  is given by Theorem  \ref{teo}. Then,  there exists a neighborhood  $V(\bar{s})$ of $\bar{s}$ in $\Lambda$ such that the following holds: every family  $\{F_t^s\}_{t \in \Sigma}$, with  $s$ in $V(\bar{s})$ verifies the above hypothesis. Theorem \ref{teo} gives us an application from  $V(\bar{s})$ to $D(t_0,K_R\eps)$
\begin{equation*}
s\mapsto \bar{t}_s
\end{equation*} 
and an application from  $V(\bar{s})$ to the set of analytical curves from the circle   $\T$ to $\C$
\begin{equation*}
s\longmapsto u^s:B_{\delta}\to \C.
\end{equation*}
The curve  $u^s$ is invariant by the fibered holomorphic dynamics  $F^s_{\bar{t}_s}$, has zero degree and its fibered rotation number equals  $\beta$. The goal of this Section is to show the
\begin{theo}\label{teopa}
The applications $s\to \bar{t}_s$, $s\to u^s$ defined above are holomorphic functions.
\end{theo}
The proof is obtained by making a more carefully study of the method (previous process and iterative algorithm), and more precisely by the ability of performing  this method in an uniform way for the whole set of families   $\{F_t^s\}_{t\in \Sigma}$, $s\in V(\bar{s})$. Let $\bar{n}$ be the unique natural number verifying 
\begin{equation*}
\frac{\delta l_{\bar{n}+1}}{24}<\eps(w_{n^*})^{-1}\leq \frac{\delta l_{\bar{n}}}{24}.
\end{equation*} 
We can apply the previous process to every family $\{F_t^s\}_{t\in \Sigma}$,    $s$ in $V(\bar{s})$. This gives us the functions  $\{(u_{0,t}^{n,i})^s\}_{(n^*\leq n <\bar{n};0\leq i\leq i_n)}$ defined for every $t$ belonging to  $D(t_0,K_R\eps)$ (here the index are uniform for every parameter $s$ in $V(\bar{s})$); this  also gives us the functions  $(\rho_{0,t}^{\bar{n}})^s,(\rho_{1,t}^{\bar{n}})^s,(\rho_t^{\bar{n}})^s$. These functions have an appropriate size  to start with the iterative algorithm at the stage  $\bar{n}$. Furthermore, these functions are defined for every $t$ in  $D(t_{\bar{n}}^s,p_{\bar{n}})$. The parameter  $t_{\bar{n}}^s$ is the only simple zero for  $\int_{\T}(\rho_{1,t}^{\bar{n}})^sd\theta$ given by the last operation in the previous procedure. Each stage  $n$ of the algorithm ($n\geq \bar{n}$) produce at some moment the functions  $\{(u_{0,t}^{n,i})^s\}_{0\leq i\leq i_n},(u_{1,t}^n)^s$ (with uniform index). These are defined for every  $t$ in  $D(t_n^s,p_n)$, where  $t_n^s$ is the only simple zero for $\int_{\T}(\rho_{1,t}^{n})^sd\theta$. The algorithm also produces  the functions $(\rho_{0,t}^{n+1})^s,(\rho_{1,t}^{n+1})^s,(\rho_t^{n+1})^s$, which have an appropriate size in order to start the stage $n+1$. These functions are defined for every  $t$ in $D(t_{n+1}^s,p_{n+1})$. 

Note that from the stage  $\bar{n}$ and so on, the parameter space  disc is the only  thing that is not uniform (the center is not uniform). We will see that all the above functions are holomorphic functions with respect to the variables  $(\theta,z,s,t)$.
\subsection{The application $\bar{t}_s$}
Note first that if one consider only one family $\{F_t\}$, the parameter  $\bar{t}$ for which one finds the invariant curve equals the only element in the intersection  $\cup_{n\geq \bar{n}}D(t_n,p_n)$. As $p_n\to 0$ when $n \to \infty$, one  deduces that $\bar{t}$ equals the limit of the sequence  $t_n$ of the centers of the discs. This suggests a strategy to show that  $\bar{t}_s$ is holomorphic. We need to show that each application $s\mapsto t_n^s$ is holomorphic. In this case, the (uniform) limit  $\bar{t}_s$ is also holomorphic, since the radius  $p_n$ are uniform.

To show that  $s\mapsto t_n^s$ is holomorphic we use a following well known fact: Let  $\{g_s:\D\to \C\}_{s\in \Lambda}$ be a family of holomorphic functions depending analytically on a parameter $s$ and verifying that each function  $g_s$ has a unique zero in $\D$ which is simple. Then this zero is a holomorphic function on the parameter    $s$ (this follows directly from the classical Implicit Function Theorem). In the next Sections we   show that the function  $\int_{\T}(\rho_{1,t}^n)^sd\theta$ is holomorphic on   $s$ and $t$.
\subsection{The functions $(u_{0,t}^{n,i})^s$,  for $n^*\leq n<\bar{n}$}
  For every  $s$ in $V(\bar{s})$ the notation $(u_{0,t}^{n,i})^s$ has an uniform sense with respect to the index, provided that $n^*\leq n< \bar{n}$. This denotes the solution of the equation of   $u_0$ for the $2^n$-order  truncate  of the function  $(\eta_{n,i}^t)^s$. Then   $(u_{0,t}^{n,i})^s$ is a trigonometric polynomial of degree  $2^n$. Its coefficients depend linearly on the Fourier series coefficients of   $(\eta_{n,i}^t)^s$. Since $(\eta_{n^*,0}^t)^s$ is holomorphic with respect to $s$ and $t$, the function $(u_{0,t}^{n^*,0})^s$ is also holomorphic with respect to $s$ and $t$. 
 
 We denote by $\tilde{<}$ the (direct) lexicographical order on $\N\times \N$. Notice that an index $(m,j)$ appears in the algorithm before another index $(n,i)$  if and only if $(m,j)\tilde{<}(n,i)$. Let  $(n,i)$ be an index. If for every index $(m,j)\tilde{<}(n,i)$ we have that  $(\eta_{m,j}^t)^s$ is holomorphic with respect to  $s$ and $t$, then $(u_{0,t}^{n,i})^s$ is holomorphic with respect to  $s$ and $t$. By an inductive reasoning we obtain  that all of the functions  $\{(u_{0,t}^{n,i})^s\}_{(n^*\leq n<\bar{n};0\leq i\leq i_n)}$, as well as the functions  $(\rho_{0,t}^{\bar{n}})^s,(\rho_{1,t}^{\bar{n}})^s,(\rho_t^{\bar{n}})^s$, are holomorphic with respect to  $s$ and $t$. In particular, the function $\int_{\T} (\rho_{1,t}^{\bar{n}})^sd\theta$ is, and consequently  $t^*_s=t^s_{\bar{n}}$, is holomorphic with respect to  $s$.
 \subsection{An arbitrary stage $n\geq \bar{n}$}
 Suppose that for every  index $(m,j)\tilde{<}(n,0)$ all of the functions $(\rho_{0,t}^{n})^s,(\rho_{1,t}^{n})^s,(\rho_t^{n})^s,\{(u_{0,t}^{m,j})^s\}$ are holomorphic on  $s$ and $t$ . Suppose also that for every  $m\leq n$ the centers of the parameter space discs  $t_m^s$ are holomorphic functions on  $s$.  For every $s$ in $V(\bar{s})$, when solving the first equation of $u_0$ in the stage  (the step $i=0$),  the trigonometrical polynomial solution  $(u_{0,t}^{n,0})^s$ also depends in an holomorphic way  on $s$ and $t$. As a consequence we have that  the functions  $(\rho_{0,t}^{n,1})^s,(\rho_{1,t}^{n,1})^s,(\rho_t^{n,1})^s$ are holomorphic.  By repeating the argument we can get the same conclusion for every function  $(u_{0,t}^{n,j})^s_{0\leq j \leq i_n},(u_{1,t}^n)^s$. This  implies that the same holds for the functions   $(\rho_{0,t}^{n+1})^s,(\rho_{1,t}^{n+1})^s,(\rho_t^{n+1})^s$. In particular, the functions 
 \begin{equation}
 \int_{\T}(\rho_{1,t}^{n+1})^sd\theta,
 \end{equation}
which are holomorphic on $t\in D(t_{n}^{\bar{n}},p_n)$, are also holomorphic on $s$.  Then there exists a holomorphic application $s\mapsto {t}_{n+1}^s$ such that
\begin{equation}
\int_{\T}(\rho_{1,{t}^s_{n+1}})^sd\theta=0.
\end{equation}
In this way, we show inductively on $n$ that the center applications   $s\mapsto t_n^s$ are holomorphic. Thus, the limit application $s\mapsto \bar{t}_s$ also corresponds to a holomorphic function on  $s$. The invariant curve $u_s$ is written as (see \ref{formuladeucero})
\begin{equation}
u^s=\sum_{n=n^*}^{\bar{n}-1}\sum_{i=0}^{i_n}(u_{0,\bar{t}_s}^{n,i})^s+\lim_{\tilde{n}\to \infty} \sum_{n=\bar{n}}^{\tilde{n}}(u_{0,\bar{t}_s}^n)^s\Big(\prod_{j=\bar{n}}^{n}(u_{1,\bar{t}_s}^j)^s\Big).
\end{equation}
The first term corresponds to the previous preparative process. This is a finite sum of holomorphic functions on $s$. The functions arising from the iterative algorithm are also holomorphic on $s$. We may then conclude that the application $s\mapsto u^s$ is holomorphic. This completes the proof of the 
Theorem \ref{teopa}$\quad_{\blacksquare}$
\bibliographystyle{plain}
\bibliography{analtesis.bib}
\end{document}